\documentclass[a4paper,11pt, twoside]{article}

\usepackage{amsfonts,amsthm,amsmath,amsfonts,latexsym,amssymb}

\usepackage{verbatim}
\usepackage{authblk} 
\usepackage{mathtools}
\usepackage{enumerate}
\usepackage{mathrsfs}
\usepackage{dsfont}
\usepackage{bbm}
\usepackage{vmargin}
\usepackage[]{appendix}
\usepackage{enumitem}
\usepackage{soul}
\usepackage{array}

\usepackage[utf8]{inputenc}
\usepackage[T1]{fontenc}

\usepackage{graphicx,upref,calligra,float}
\usepackage[rflt]{floatflt}
\usepackage[svgnames]{xcolor}

\usepackage[colorlinks=true,linkcolor=azul,citecolor=azul]{hyperref}
\usepackage{fancyhdr}
\usepackage{lastpage}
\usepackage{etoolbox}
\usepackage{tocloft}

\newtheorem{teor}{Theorem}[section]

\newtheorem{teo}[teor]{Theorem}
\newtheorem*{teo*}{Theorem}
\newtheorem{lem}[teor]{Lemma}

\newtheorem{cor}[teor]{Corollary}

\theoremstyle{definition}

\newtheorem{obs}[teor]{Remark}
\newtheorem*{obs*}{Remark}
\newtheorem{prop}[teor]{Proposition}

\newtheorem{exa}[teor]{Example}

\definecolor{azul}{rgb}{0.1,0.6,0.86}
\definecolor{naranja}{RGB}{249,153,96}

\def\bdem{\begin{proof}}
\def\edem{\end{proof}}

\def\fii{\varphi }

\def\a{\alpha}

\def\N{\mathbb{N}}
\def\Z{\mathbb{Z}}
\def\R{\mathbb{R}}
\def\Q{\mathbb{Q}}
\def\C{\mathbb{C}}

\def\cI{\mathcal{I}}

\newcommand{\e}{\mathbb{E}}

\newcommand{\efe}{\mathcal{F}}
\newcommand{\pe}{\mathbb{P}}

\newcommand{\peso}[1]{ \quad \mbox{  #1 } \quad }

\newcommand{\schw}{\mathcal{S}(\R^d)}

\def\ccd{C_c(\R^d)}

\title{Diffraction of randomly perturbed atomic measures}
\author[1,2]{Jorge Antezana\thanks{Corresponding author. E-mail \texttt{jorge.antezana@ub.edu}. Supported by Grant PID2024-160033NB-I00 funded by MICIU/AEI/10.13039/501100011033 and by FEDER, UE.}}
\author[3,4]{Luciano Scazzola \thanks{E-mail \texttt{lscazzola@mate.unlp.edu.ar}. Supported by National University of La Plata Grant 11X974, CONICET-PIP 11220210100954CO. Part of this work was carried out during L. Scazzola's research stay at the Universitat de Barcelona in 2025, with airfare support from the Asociaci\'on Universitaria Iberoamericana de Postgrado (AUIP) through its Programa Magallanes de Movilidad Acad\'emica Internacional (grant PMG1-2025-12-056). }}

\affil[1]{Departament de Matem\`atiques i Inform\`atica, Universitat de Barcelona, Barcelona, Spain.}
\affil[2]{Centre de Recerca Matem\`atica, Barcelona, Spain.}
\affil[3]{Universidad Nacional de La Plata, La Plata, Argentina.}
\affil[4]{Instituto Argentino de Matem\'atica ``Alberto P. Calder\'on'', CONICET, Argentina.}

\date{}

\hypersetup{pdftitle={Diffraction of randomly perturbed atomic measures},pdfauthor={Jorge Antezana and Luciano Scazzola}}

\begin{document}

\maketitle

\begin{abstract}
We study independent, identically distributed displacements of a
weighted atomic configuration on $\R^d$. We assume that the measure 
associated to the configuration is translation bounded, has complex weights
and locally finite support, and that its autocorrelation exists. For any common displacement law other than a point mass, we prove that, with probability one, the displaced measure has autocorrelation if and only if the squared moduli of the original weights have a limiting sum per unit volume. No moment assumption on the displacements is required. This density exists automatically for uniformly discrete supports. In general, we show that even a positive measure with the same autocorrelation as the integer lattice can lose its autocorrelation under arbitrarily small independent displacements almost surely. When the density exists, we give an explicit formula for the diffraction and prove that the observed scattering measures converge to it almost surely in the sense of tempered distributions as the observation radius grows, without passing to a subsequence. For a known non-atomic displacement law with a characteristic function without zeros, a correction computed from the observed configuration recovers the original diffraction even when the ordinary limit fails. The existence criterion and diffraction formulas extend to
independent, identically distributed finite complex marks with uniformly bounded
total variation. Observable recovery also extends to atomic
marks under suitable additional assumptions.
\end{abstract}

\newpage

\tableofcontents

\newpage

\section{Introduction}

Diffraction probes long-range order through the spatial correlations of a configuration of scatterers. A measure on $\R^d$ records their positions and scattering amplitudes. Its autocorrelation describes the volume-averaged two-point statistics, and the Fourier transform of the autocorrelation gives the diffraction measure. This description applies to periodic crystals, quasicrystals and more general weighted structures \cite{Hof,BaakeGrimmBook,BaakeGrimm2012}.

Can two configurations with exactly the same diffraction respond so differently to independent displacements that only one retains the existence of autocorrelation? The answer is yes. By splitting selected atoms of the integer lattice, we obtain a configuration with exactly the same autocorrelation and hence the same diffraction. Yet under independent displacements with any common law other than a point mass, the perturbed lattice almost surely has  autocorrelation, while the perturbed split configuration almost surely has not. This holds even if the displacement law is supported in an arbitrarily small ball centered at the origin. 

Hof~\cite{HofHighTemp} proved an almost-sure diffraction formula for uniformly discrete point sets satisfying a chain-frequency condition. Baake, Birkner and Moody~\cite{BaakeBirknerMoody2010} treated independent decorations of fixed point sets of finite local complexity (FLC), allowing passage to a subsequence of observation windows. However,  for a weighted configuration, the example above shows why existence itself needs a criterion. We prove that the decisive quantity is the sum of squared moduli of the atomic weights per unit volume. We then determine what can be recovered from a single displaced configuration, including when its autocorrelation fails to exist.

More precisely, let $X\subseteq\R^d$ be a locally finite set and consider a purely atomic measure
\[
    \mu=\sum_{p\in X}\mu(p)\,\delta_p ,\qquad \mu(p)\in\C\setminus\{0\},
\]
which we assume to be translation bounded (TB), meaning that
\[
    \sup_{x\in\R^d}|\mu|(B_1(x))<\infty.
\]
Here $B_R(x)$ is the closed ball of radius $R$ centered at $x$, and $B_R:=B_R(0)$. The weights are complex scattering amplitudes. Translation boundedness controls the total absolute weight in each unit ball and allows the number of atoms in such balls to be unbounded. Write $m_d$ for Lebesgue measure and $\mu_R:=\mu|_{B_R}$. The finite-volume autocorrelation is
\[
    \gamma_\mu^{(R)}:=\frac{\mu_R*\widetilde{\mu_R}}{m_d(B_R)}.
\]
Its vague limit as $R\to\infty$ through all radii, when it exists, is the \emph{natural autocorrelation} $\gamma_\mu$. We use this averaging convention throughout. The atoms of the diffraction measure $\widehat{\gamma_\mu}$ are the Bragg peaks, and its continuous part describes diffuse scattering. Section~\ref{prelim} gives the definitions. We assume that $\gamma_\mu$ exists when studying existence and recovery after displacement. The local estimates and the uniform bound of Theorem~\ref{thm:uniform-TB} also apply to backgrounds for which this limit fails.

\subsection{Main results}

We begin with the pure displacement model, in which only the positions of the atoms are perturbed. Let
$\{\xi_p\}_{p\in X}$ be an i.i.d.\ family of $\R^d$-valued random vectors
with common law $\pe_\xi$, defined on a probability space
$(\Omega,\efe,\pe)$. No assumption is made on $\pe_\xi$; in particular, no moment of $\xi$ is required to be finite. The perturbed random measure is
\[
    \mu_\xi
    :=
    \sum_{p\in X}\mu(p)\,\delta_{p+\xi_p} \,.
\]
Throughout, an almost-sure statement about vague convergence means that, on a single event of probability one, the convergence holds for every $f\in\ccd$.

\paragraph{Existence of the displaced autocorrelation.}
For each $R$, a pair of sites $p,q\in X\cap B_R$ contributes to the mass of $\gamma_\mu^{(R)}$ at the origin precisely when $p=q$. Thus
\[
    \gamma_\mu^{(R)}(\{0\})
    =\frac{\boldsymbol E_R}{m_d(B_R)},
    \qquad
    \boldsymbol E_R:=\sum_{p\in X\cap B_R}|\mu(p)|^2.
\]
We call $\boldsymbol E_R$ the \emph{diagonal energy} in $B_R$ and denote the limit of $\boldsymbol E_R/m_d(B_R)$, when it exists, by $\boldsymbol E$, the \emph{diagonal energy density}. Vague convergence of $\gamma_\mu^{(R)}$ alone does not imply convergence of its mass at a singleton. If $X$ is uniformly discrete, a single continuous test function supported in a sufficiently small ball and equal to one at the origin isolates the diagonal at every radius. In that case $\boldsymbol E$ exists and equals $\gamma_\mu(\{0\})$. Without uniform separation, differences between distinct sites can approach zero and contribute to the atom of the limiting autocorrelation.

Independent displacements act differently on diagonal and off-diagonal pairs. Each self-pair still has difference zero after displacement. For $p\ne q$, the original difference $p-q$ is shifted by $\xi_p-\xi_q$, whose law is $\pe_\xi*\widetilde{\pe_\xi}$. This distinction leads to the following criterion.

\begin{teo}\label{autocoraleatoriaexists}
Suppose that $\gamma_\mu$ exists.
\begin{enumerate}[label=\textup{(\roman*)}]
\item If the diagonal energy density $\boldsymbol{E}$ exists, then almost surely the autocorrelation of $\mu_\xi$ exists and
\begin{equation}\label{eq:random-auto-intro}
    \gamma_{\mu_\xi}
    =\gamma_\mu*(\pe_\xi*\widetilde{\pe_\xi})
    +\boldsymbol{E}\,\bigl(\delta_0-\pe_\xi*\widetilde{\pe_\xi}\bigr).
\end{equation}
\item If $\boldsymbol{E}$ does not exist and $\pe_\xi$ is not a point mass, then almost surely the finite-volume autocorrelations of $\mu_\xi$ have no vague limit as $R\to\infty$ through all radii.
\end{enumerate}
\end{teo}

For a point-mass law the perturbation is a translation, so $\gamma_{\mu_\xi}=\gamma_\mu$ even when $\boldsymbol E$ fails to exist.
To see the obstruction in Theorem~\ref{autocoraleatoriaexists}, return to the integer lattice. Choose a subset $A\subseteq\Z$ for which $\#(A\cap[-R,R])/(2R)$ has no limit and set
\[
    \mu=\sum_{n\in\Z\setminus A}\delta_n
    +\frac12\sum_{n\in A}
      \bigl(\delta_n+\delta_{n+\varepsilon_n}\bigr),
    \qquad \varepsilon_n=\frac1{|n|+2}.
\]
At a split site, the two self-pairs have total weight $1/2$ at zero. The two cross-pairs contribute the remaining weight, $1/4$ at each of $\pm\varepsilon_n$. As $|n|\to\infty$, a continuous test function sees the same total contribution as from the original atom. Applying the same reasoning to pairs of cells with a fixed integer separation gives $\gamma_\mu=\delta_{\Z}$ (Example~\ref{ex:diagonal-energy}). Each split reduces the diagonal energy by $1/2$. Accounting for the bounded contribution from the window boundary gives
\[
    \frac{\boldsymbol E_R}{2R}
    =1-\frac{\#(A\cap[-R,R])}{4R}+o(1).
\]
The right-hand side oscillates. The unsplit lattice has diagonal energy density one. The criterion thus gives opposite answers for two positive configurations with identical diffraction and at most two sites per unit cell $[n,n+1)$.

The criterion applies to many configurations that are not uniformly discrete. For instance, if $a>0$ is irrational, then $\delta_{\Z}+\delta_{a\Z}$ has autocorrelation and its diagonal energy density is $1+a^{-1}$. 
Signed weights are allowed as well: the Thue--Morse and Rudin--Shapiro combs, with weights in $\{-1,1\}$ on $\Z$, have $\boldsymbol E=1$; see \cite[Secs.~5.1--5.2]{BaakeGrimm2012}. The Rudin--Shapiro comb is a striking case, since its autocorrelation is $\delta_0$: by \eqref{eq:random-auto-intro}, $\gamma_{\mu_\xi}=\delta_0$ for every displacement law, so its diffraction is not affected at all by independent displacements.

\paragraph{Observed diffraction and corrected limits.}

The scattering intensity in a bounded observation window is determined by the atoms at their final positions. We therefore form
\begin{equation}\label{eq:gamma}
    \gamma_{\mu_\xi}^{(R)}
    :=\frac{(\mu_\xi)_R*\widetilde{(\mu_\xi)_R}}{m_d(B_R)},
    \qquad R\geq1.
\end{equation}
Their Fourier transforms are positive scattering measures, with densities
\[
    \frac{1}{m_d(B_R)}
    \bigl|\widehat{(\mu_\xi)_R}(\lambda)\bigr|^2
\]
with respect to $m_d$. Restricting to the window and displacing do not commute. Displacing the atoms that were originally in $B_R$ gives
\[
    (\mu_R)_\xi=\sum_{p\in X\cap B_R}\mu(p)\,\delta_{p+\xi_p},
\]
whereas $(\mu_\xi)_R$ consists of the displaced atoms that lie in $B_R$, wherever they came from. Only $(\mu_\xi)_R$ can be computed from the final, unlabeled configuration, while $(\mu_R)_\xi$ is the object to which laws of large numbers apply directly. Corollary~\ref{primersoportes} shows that, for every displacement law, the difference between their normalized two-point measures tends vaguely to zero, even when neither family converges:
\begin{equation}\label{eq:observed-window-intro}
    \frac{(\mu_\xi)_R*\widetilde{(\mu_\xi)_R}
          -(\mu_R)_\xi*\widetilde{(\mu_R)_\xi}}{m_d(B_R)}
    \xrightarrow[R\to\infty]{\mathrm{vague}}0
    \qquad\text{almost surely}.
\end{equation}

The next result describes what converges when the diagonal energy density fails to exist. The correction term involves the diagonal energy of the atoms observed in the window,
\[
 \boldsymbol{E}_{\xi,R}
    :=
    \sum_{p\in X:\,p+\xi_p\in B_R}|\mu(p)|^2 \,,
\]
which, unlike $\boldsymbol E_R$, can be read off from the displaced configuration whenever no two displaced atoms coincide: it is then equal to $\sum_{x\in B_R}|\mu_\xi(\{x\})|^2$.

\begin{teo}\label{autocoraleatoria}
Suppose that $\gamma_\mu$ exists. Then, almost surely,
\begin{equation}\label{autocorlim}
    \frac{
        (\mu_\xi)_R*\widetilde{(\mu_\xi)_R}
        -\boldsymbol{E}_{\xi,R}\,
        \bigl(\delta_0-\pe_\xi*\widetilde{\pe_\xi}\bigr)
    }{m_d(B_R)}
    \xrightarrow[R\to\infty]{\mathrm{vague}}
    \gamma_\mu*(\pe_\xi*\widetilde{\pe_\xi}).
\end{equation}
\end{teo}

The boundary estimates give $(\boldsymbol E_{\xi,R}-\boldsymbol E_R)/m_d(B_R)\to0$ almost surely. Thus $\boldsymbol E_R$ can replace $\boldsymbol E_{\xi,R}$ in~\eqref{autocorlim}. Along any sequence for which the energy per unit volume tends to $e$, the finite-volume autocorrelations therefore converge almost surely to
$\gamma_\mu*(\pe_\xi*\widetilde{\pe_\xi})+e(\delta_0-\pe_\xi*\widetilde{\pe_\xi})$.
Every value between the lower and upper limits of $\boldsymbol E_R/m_d(B_R)$ occurs in this way (Remark~\ref{rem:accumulation-segment}).

Passing from the observed autocorrelations to scattering intensities requires control at large distances as well as vague convergence. The displaced measure can lose translation boundedness even when the original configuration is a lattice (Example~\ref{ex:lattice-gaussian-not-TB}). Its finite-volume autocorrelations nevertheless satisfy the following uniform bound.

\begin{teo}\label{thm:uniform-TB}
Let $\mu$ be a translation bounded, purely atomic measure with locally finite support, and let the displacement law be arbitrary. Then, almost surely, the family $\{\gamma_{\mu_\xi}^{(R)}:R\geq1\}$ is uniformly translation bounded: for almost every $\omega\in \Omega$ there is a constant $C(\omega)<\infty$ such that
\begin{equation}\label{eq:main-TB}
    \sup_{R\geq1}\,\sup_{x\in\R^d}
    |\gamma_{\mu_\xi}^{(R)}|(B_1(x))\leq C(\omega).
\end{equation}
\end{teo}

When $\gamma_\mu$ and $\boldsymbol E$ exist, this bound gives convergence in $\mathcal S'(\R^d)$ (Lemma~\ref{lem:vague-to-distribution}), where the Fourier transform is continuous. The diffraction can thus be obtained from observations of a single realization as the window grows. Write $\varphi_\xi(\lambda):=\e(e^{-2\pi i\langle\lambda,\xi\rangle})$ for the characteristic function of the displacement law.

\begin{teo}\label{formulasautoanddifrac}
Assume that $\gamma_\mu$ and $\boldsymbol E$ exist. Then, almost surely, the diffraction measure of $\mu_\xi$ is
\begin{equation}\label{eq:random-diff-intro}
    \widehat{\gamma_{\mu_\xi}}
    =|\varphi_\xi|^2\,\widehat{\gamma_\mu}
    +\boldsymbol{E}\,\bigl(1-|\varphi_\xi|^2\bigr)\,m_d ,
\end{equation}
and, for every $f\in\schw$,
\begin{equation}\label{eq:computable-diffraction}
    \lim_{R\to\infty} \widehat{\gamma_{\mu_\xi}^{(R)}}(f)=\widehat{\gamma_{\mu_\xi}}(f)  \,.
\end{equation}
\end{teo}

For Gaussian displacements with covariance matrix $\Sigma$, the multiplier $|\varphi_\xi(\lambda)|^2$ is the Debye--Waller intensity factor $\exp(-4\pi^2\lambda^{\mathsf T}\Sigma\lambda)$. It attenuates the original diffraction, including any continuous component already present. The additional scattering generated by the displacements is absolutely continuous, with density $\boldsymbol E(1-|\varphi_\xi|^2)$. In the unweighted, uniformly discrete case $\boldsymbol E$ is the point density, and~\eqref{eq:random-diff-intro} recovers Hof's formula~\cite{HofHighTemp}.

\paragraph{Recovery from a single realization.}
A Bragg peak at $\lambda$ retains the fraction $|\varphi_\xi(\lambda)|^2$ of its intensity. The absolutely continuous correction creates no new peaks. This gives a direct recovery formula for pure point diffraction.

\begin{cor}\label{purelyatomicrecovery}
Assume that $\gamma_\mu$ and $\boldsymbol{E}$ exist, that $\widehat{\gamma_\mu}$ is a pure point measure, and that $\varphi_\xi$ has no zeros. Then, almost surely,
\[
  \widehat{\gamma_\mu}
  =|\varphi_\xi|^{-2}\,
      (\widehat{\gamma_{\mu_\xi}})_{\mathrm{pp}}.
\]
Consequently, if the displacement law is known, the original diffraction is determined by a single realization of $\mu_\xi$, through the limit~\eqref{eq:computable-diffraction}.
\end{cor}

For \emph{arbitrary} original diffraction there is a second route, which works even when $\boldsymbol E$ fails to exist. If the displacement law is non-atomic, almost surely no two displaced atoms coincide. On that event the diagonal energy $\boldsymbol E_{\xi,R}$ is exactly the sum of the squared moduli of the observed weights in $B_R$, simultaneously for every $R$. Corollary~\ref{cor:observable-recovery} then gives, almost surely,
\[
  \widehat{\gamma_{\mu_\xi}^{(R)}}
  -\frac{\boldsymbol{E}_{\xi,R}}{m_d(B_R)}\,
   \bigl(1-|\varphi_\xi|^2\bigr)\,m_d
  \xrightarrow[R\to\infty]{\mathcal S'(\R^d)}
  |\varphi_\xi|^2\,\widehat{\gamma_\mu}.
\]
If the law is known and $\varphi_\xi$ is nowhere zero, division of the limiting measure by $|\varphi_\xi|^2$ on compact sets recovers $\widehat{\gamma_\mu}$. For example, a Cauchy displacement law on $\R$ has no first moment, has no atoms and has a characteristic function with no zeros. The correction therefore recovers the common diffraction of the integer lattice and its split version even though only the former has an ordinary diffraction limit after displacement (Section~\ref{randomautocorrelation}).

\paragraph{Independent atomic marks.}
Pure displacements were treated first because they show the mechanism with the least notation. Random occupation, random amplitudes and the splitting of an atom into several weighted atoms require a more general perturbation, in which the atom at $p$ is replaced by an independent random atomic measure $\Theta_p$, translated to $p$. We call $p$ a \emph{parent}, $\Theta_p$ its \emph{mark}, and the individual weighted atoms in its defining sum
its \emph{descendants}; this is the fixed-background analogue of a cluster point process \cite[Sec.~6.3]{DaleyVereJonesI}.

For each $p\in X$, let
\[
\Theta_p:=
\sum_{k=1}^{N(p)}
\alpha_{p,k}\,\delta_{\xi_{p,k}},
\qquad N(p)\in\mathbb N\cup\{\infty\},\quad \alpha_{p,k}\in \C,\quad \xi_{p,k}\in\R^d .
\]
We assume that the families
\[
\bigl(N(p),
(\alpha_{p,k},\xi_{p,k})_{1\leq k\leq N(p)}\bigr),
\qquad p\in X,
\]
are i.i.d.; the number, the weights and the positions of the descendants of a single parent may depend on one another in an arbitrary way. Let $\Theta=\sum_{k=1}^{N}\alpha_k\,\delta_{\xi_k}$ denote a generic mark. Our only assumption is that the total absolute descendant weight is uniformly bounded: there is a constant $C<\infty$ such that, almost surely,
\begin{equation}\label{eq:mark-bound}
    \sum_{k=1}^{N}|\alpha_{k}|\leq C \,.
\end{equation}
Thus each mark is a finite complex measure, even when $N=\infty$ and its support has accumulation points, and no moment of the positions $\xi_k$ is required. As shown at the end of Section~\ref{generalcase},
the atomic structure of the marks is not essential
for the existence criterion or the autocorrelation
and diffraction formulas: these results extend to
general finite complex marks with uniformly bounded
total variation. This boundedness assumption is stronger
than the square-integrability condition on the total
variation in \cite{BaakeBirknerMoody2010}, while our results
allow more general fixed backgrounds and establish
almost-sure convergence along all radii. The atomic setting is nevertheless particularly relevant
to the recovery problem, since under suitable conditions
the correction can be computed directly from the observed
atomic weights.

Let $\tau_p$ denote translation by $p$. The perturbed measure is
\[
    \mu_\Theta
    :=\sum_{p\in X}\mu(p)\,\tau_p\Theta_p
    =\sum_{p\in X}\sum_{k=1}^{N(p)}\mu(p)\,\alpha_{p,k}\,\delta_{p+\xi_{p,k}},
\]
and its finite-volume autocorrelations are $\gamma_{\mu_\Theta}^{(R)}:=(\mu_\Theta)_R*\widetilde{(\mu_\Theta)_R}/m_d(B_R)$, $R\geq1$. Two finite measures govern the result,
\[
 M:=\e(\Theta*\widetilde\Theta)
 \qquad\text{and}\qquad
 \kappa:=\e(\Theta)*\widetilde{\e(\Theta)} \,:
\]
$M$ records the correlations between descendants of the same parent, and $\kappa$ those between descendants of different parents. Their Fourier transforms are $\widehat M(\lambda)=\e\bigl(|\widehat\Theta(\lambda)|^2\bigr)$ and $\widehat\kappa(\lambda)=|\fii_\Theta(\lambda)|^2$, where
\[
    \fii_\Theta(\lambda):=\e\bigl(\widehat\Theta(\lambda)\bigr)
    =\e\Bigl(\sum_{k=1}^{N}\alpha_k \,e^{-2\pi i\langle\lambda,\xi_k\rangle}\Bigr) ,
\]
so that $M-\kappa$ has Fourier transform
\[
    \operatorname{Var}(\widehat{\Theta})(\lambda)
    :=\e\bigl(|\widehat{\Theta}(\lambda)|^2\bigr)
      -|\fii_\Theta(\lambda)|^2 ,
\]
the variance of the random variable $\widehat\Theta(\lambda)$. Finally, let
\[
 c_\Theta:=\e\Bigl(\sum_{k=1}^{N} |\alpha_{k}|^2\Bigr) \qquad \text{and} \qquad \boldsymbol{E}_{\Theta,R} := \frac{1}{c_\Theta} \sum_{p\in X} |\mu(p)|^2 \sum_{k=1}^{N(p)} |\alpha_{p,k}|^2\, \mathbbm{1}_{\{p+\xi_{p,k}\in B_R\}}\,,
\]
where we assume $c_\Theta>0$, since otherwise $\mu_\Theta=0$ almost surely. The normalization ensures that
$\boldsymbol{E}_{\Theta,R}-\boldsymbol{E}_R
=o(m_d(B_R))$ almost surely
(Lemma~\ref{energyclustergeneral}). We use the notation $\gamma_{\mu_\Theta}^{(R)}$ of~\eqref{eq:gamma} for the marked observation window.

\begin{teo}\label{teo:cluster-main-intro}
Assume that $\gamma_\mu$ exists and that the i.i.d.\ marks satisfy \eqref{eq:mark-bound}, with $c_\Theta>0$.
\begin{enumerate}[label=\textup{(\roman*)}]
\item Almost surely,
\begin{equation}\label{firstpart}
    \frac{
        (\mu_\Theta)_R*\widetilde{(\mu_\Theta)_R}
        -\boldsymbol{E}_{\Theta,R}\,
        \bigl(M-\kappa\bigr)
    }{m_d(B_R)}
    \xrightarrow[R\to\infty]{\mathrm{vague}}
    \gamma_\mu*\kappa.
\end{equation}
\item If the diagonal energy density $\boldsymbol{E}$ of $\mu$ exists, then almost surely the autocorrelation of $\mu_\Theta$ exists and
\begin{equation}\label{secondpart}
    \gamma_{\mu_\Theta}=\gamma_\mu*\kappa+\boldsymbol{E}\,(M-\kappa),
\end{equation}
so that its diffraction measure is
\begin{equation}\label{thirdpart}
    \widehat{\gamma_{\mu_\Theta}}
    =|\fii_\Theta|^2\,\widehat{\gamma_\mu}
    +\boldsymbol{E}\,\operatorname{Var}(\widehat{\Theta})\,m_d .
\end{equation}
\item In the situation of \textup{(ii)}, almost surely, for every $f\in\schw$,
\begin{equation}\label{eq:computable-cluster-diffraction}
    \lim_{R\to\infty}\widehat{\gamma_{\mu_\Theta}^{(R)}}(f)
    =\widehat{\gamma_{\mu_\Theta}}(f).
\end{equation}
\end{enumerate}
\end{teo}

The pure displacement model is the case $\Theta=\delta_\xi$, for which $c_\Theta=1$ and
\[
    M=\delta_0,\qquad \kappa=\pe_\xi \ast \widetilde{\pe_\xi}, \qquad
    \fii_\Theta=\varphi_\xi,\qquad
    \operatorname{Var}(\widehat{\Theta})=1-|\varphi_\xi|^2 ,
\]
which recovers the displacement formulas. The choice $\Theta=\alpha\,\delta_0$, with $\alpha$ a bounded random variable, gives random amplitudes. When $\alpha$ is Bernoulli with mean $\theta$, each atom is kept with probability $\theta$ and \eqref{thirdpart} becomes
\[
    \widehat{\gamma_{\mu_\Theta}}=\theta^2\,\widehat{\gamma_\mu}+\boldsymbol{E}\,\theta(1-\theta)\,m_d ,
\]
the classical formula for the random occupation model; compare \cite[Remark~10]{BaakeBirknerMoody2010}.

The same energy obstruction persists for marks. If $\Theta$ is not almost surely equal to a fixed measure and $\boldsymbol E$ does not exist, then almost surely $\mu_\Theta$ has no autocorrelation (Corollary~\ref{cor:marked-obstruction}). As for displacements, the corrected limit \eqref{firstpart} recovers $\widehat{\gamma_\mu}$ from a single realization when the law of the marks is known and $\fii_\Theta$ has no zeros, provided that $\boldsymbol E_{\Theta,R}$ can be computed from the observed measure. This is the case when almost surely no two descendants occupy the same position, since then $c_\Theta\,\boldsymbol E_{\Theta,R}=\sum_{x\in B_R}|\mu_\Theta(\{x\})|^2$. When descendants may coincide, their weights are added before the modulus is squared, and the observed weights need not determine the correction. Section~\ref{generalcase} discusses this observable
recovery and shows that, for general finite marks,
the corrected limit remains valid with the original
diagonal energy $\boldsymbol E_R$ in place of
$\boldsymbol E_{\Theta,R}$.

\paragraph{Local regularity.}
Since independent displacements can destroy translation boundedness even for a lattice (Example~\ref{ex:lattice-gaussian-not-TB}), it is natural to ask how large the local mass of $\mu_\xi$ can become. Let $Q=[-1/2,1/2)^d$. Theorem~\ref{emeR} shows that, for every displacement law and every $D>d\,\sup_{p\in X}|\mu(p)|$, almost surely
\begin{equation}\label{eq:quasi-TB-intro}
    |\mu_\xi|(Q+q)
    \leq D\,\frac{\log\|q\|}{\log\log\|q\|}
    \qquad\text{for all $q\in\Z^d$ with $\|q\|$ sufficiently large,}
\end{equation}
where the threshold is random. In particular, $\mu_\xi$ is almost surely quasi-translation bounded in the sense of Section~\ref{randomperturbTB}. The order of growth is sharp along a sequence of unit cubes for the lattice $\Z^d$ with a displacement law that has polynomial tails and finite first moment (Proposition~\ref{ex:sharp-local-growth}). The estimate does not require $\gamma_\mu$ to exist, and it extends to marks satisfying~\eqref{eq:mark-bound} (Theorem~\ref{emeRgeneral}).

\subsection{Relation with previous results}\label{subsec:previous-work}

Table~\ref{tab:comparison} summarizes the hypotheses and modes of convergence for fixed backgrounds.

\begin{table}[ht]
\centering
\small
\renewcommand{\arraystretch}{1.25}
\begin{tabular}{@{}>{\raggedright\arraybackslash}p{2.3cm}>{\raggedright\arraybackslash}p{3.6cm}>{\raggedright\arraybackslash}p{4.2cm}>{\raggedright\arraybackslash}p{4.4cm}@{}}
\hline
 & Hof~\cite{HofHighTemp} & Baake, Birkner and Moody~\cite[Thm.~2 and Rem.~11]{BaakeBirknerMoody2010}& This paper\\
\hline
Background & unweighted, uniformly discrete, chain frequencies & 
unweighted FLC point set with autocorrelation;
more generally, pure point autocorrelation with
locally finite support
 & translation bounded, complex weights, locally finite support, with autocorrelation\\
Perturbation & i.i.d.\ displacements with finite first moment & i.i.d.\ finite random measures with square-integrable total variation & i.i.d.\ displacements with arbitrary law; i.i.d.\ atomic marks satisfying~\eqref{eq:mark-bound}; \\
Window in the proof & original positions & compared with the observed window in mean & compared with the observed window almost surely\\
Convergence & almost sure, vague & almost sure, possibly on a subsequence, vague & almost sure along all radii, vague and in
$\mathcal S'(\R^d)$\\
\hline
\end{tabular}
\caption{Almost-sure diffraction results for independent perturbations of a fixed background.}
\label{tab:comparison}
\end{table}

Hof~\cite[Secs.~2--3]{HofHighTemp} obtained formula~\eqref{eq:random-diff-intro} for unweighted, uniformly discrete point sets and displacements with finite first moment, under the following hypothesis, which we call \emph{Hof's chain-frequency condition}. For every translate $A$ of a sufficiently small cube, the pairs $(x,y)\in X\times X$ with $x-y\in A$ are arranged in chains, and the pairs of odd rank and those of even rank in their chain are required to have separate frequencies. The purpose of this condition is to split a sum of dependent terms into two sums of independent ones; for several types of atoms, a strengthened version is needed \cite[Sec.~4]{HofHighTemp}. In the argument of Baake, Birkner and Moody each class whose size tends to infinity has the same limiting mean per pair, so only the total pair frequency enters \cite[Eqs.~(31)--(34)]{BaakeBirknerMoody2010}. In Hof's setting and in the FLC setting of
\cite[Theorem~2]{BaakeBirknerMoody2010}, the support is
uniformly discrete, so the diagonal energy density is
the point density $\gamma_\mu(\{0\})$, and the phenomenon
described in Theorem~\ref{autocoraleatoriaexists}
cannot occur.

One point in this comparison concerns the observation window. Hof defines the autocorrelation of the displaced configuration through its restriction to a cube $C_L$, but the sums in the proof of his formula~(3.1) run over the pairs $x,y\in X\cap C_L$, that is, over atoms selected by their original positions \cite[pp.~59--60]{HofHighTemp}. Since atoms may enter or leave the cube, the two families of approximants differ, and identifying their limits requires a bound on the pairs that involve a crossing atom. Corollary~\ref{primersoportes} provides this bound for arbitrary displacement laws; see~\eqref{eq:observed-window-intro}. Baake, Birkner and Moody address the same issue in their Proposition~2, which compares the observed pair sum in mean with a sum in which one factor is selected by the parent sites. This estimate uses no spatial moment and permits a passage to a subsequence for almost-sure convergence. Their Remark~11 indicates an extension to backgrounds whose autocorrelation is pure point with locally finite support, and their Sections~5.3--5.4 obtain almost-sure limits by a different route, for backgrounds sampled from a stationary ergodic point process. In the present paper the background is a fixed measure, the comparison holds almost surely along all radii, and it remains valid when the diagonal energy oscillates. Moreover, while Baake, Birkner and Moody assume FLC
in their Theorem~2, Hof's result does not require
FLC or pure point diffraction, but instead relies
on the chain-frequency condition described above.
Our results require neither FLC nor this
chain-frequency condition. For instance, the
sine-modulated lattice of Example~\ref{ex:sine-modulation}
is uniformly discrete but not FLC, and its
autocorrelation has absolutely continuous components.

The coherent and diffuse terms in \eqref{eq:random-diff-intro} are familiar from structure-factor calculations for randomly displaced lattices and for splitting processes~\cite{Gabrielli,KlattKimTorquato,GabrielliJoyce2008}. K\"ulske~\cite{Kuelske2003,KuelskeCMP2003} obtained finite-volume concentration bounds for separated backgrounds, including models with dependent disorder; they are complementary to the almost-sure limits along all radii obtained here. Petrache and Viera~\cite[Sec.~5.2]{PetracheViera2023} ask under what conditions an almost-sure formula relates the diffraction of a randomly displaced point set to its original diffraction. For i.i.d.\ displacements, Theorem~\ref{formulasautoanddifrac} answers this question in the sense of tempered distributions. Related comparisons for linear observables appear in~\cite{Yakir2022,PetracheViera2023}, and the general theory of reflected Eberlein convolutions and diffraction is developed in~\cite{LenzSpindelerStrungaru2024,LenzStrungaru2026}.

Our earlier result for crystalline measures uses spatial averages of complex amplitudes to determine $\widehat\mu$ wherever the known displacement factor is nonzero, and hence to recover $\mu$ when that factor is nowhere zero \cite[Theorem~1.2 and Corollary~1.3]{nosotros}. Here the observed quantities are quadratic. For a weighted background with autocorrelation, the diagonal energy governs existence of the perturbed limit, and the corrected intensities yield the original diffraction under the hypotheses of Corollary~\ref{cor:observable-recovery}. The sharp local mass estimate also strengthens the global growth control in~\cite[Remark~4.3]{nosotros}. It applies to arbitrary weighted TB backgrounds, complementing local estimates for specific random configurations such as~\cite[Lemma~2.2]{GerminetHislopKlein2006}.

Finally, hyperuniformity can be read off from the behavior of the diffraction measure near the origin \cite{BaakeGrimmHU,BjorklundHartnickHU,BjorklundBylehn2026}, and the effect of random perturbations on the hyperuniformity of lattices, including dependent displacements, was studied in \cite{DereudreFlimmelHuesmannLeble}. When the diagonal energy density exists, our formulas
show that independent displacements preserve spectral
hyperuniformity for every background considered here.
More generally, if the original measure is spectrally
hyperuniform and $\boldsymbol E>0$, independent finite
marks with uniformly bounded total variation preserve
spectral hyperuniformity if and only if the total mass
$\Theta(\R^d)$ of a mark is almost surely constant
(Corollary~\ref{cor:spectral-hu} and its extension
at the end of Section~\ref{generalcase}).

\subsection{Organization of the paper}

Section~\ref{prelim} collects the notation and the basic facts on
translation bounded measures, vague and distributional convergence,
autocorrelation and diffraction, including the criterion from
\cite{nosotros} that upgrades vague convergence to convergence in
$\mathcal{S}'(\R^d)$, as well as some probabilistic preliminaries.
Section~\ref{randomperturbTB} establishes local mass and boundary-crossing
estimates, some of them in a general form that also covers marks.
Section~\ref{averaging} proves the window comparison and the strong law
for pair sums in a common framework for displacements and marks.
Sections~\ref{randomautocorrelation} and~\ref{randomdiffraction}
prove Theorems~\ref{autocoraleatoriaexists}--\ref{formulasautoanddifrac}
and the recovery results, including a general uniform translation
boundedness estimate for normalized random autocorrelations.
Section~\ref{generalcase} extends the results, including
observable recovery under suitable hypotheses, to
independent atomic marks, discusses the consequences
for spectral hyperuniformity, and finally extends the
existence criterion, autocorrelation and diffraction
formulas to general finite marks.

\section{Preliminaries}\label{prelim}

Throughout the paper, $\|\cdot\|$ denotes the Euclidean norm on $\R^d$ and
$\langle\cdot,\cdot\rangle$ its standard inner product. For $x_0\in\R^d$ and
$R>0$, we write
\[
    B_R(x_0):=\{x\in\R^d:\|x-x_0\|\leq R\},
    \qquad
    B_R:=B_R(0),
\]
and $m_d$ denotes Lebesgue measure on $\R^d$. If $Y\subseteq\R^d$, we use the
shorthand
\[
    Y_R:=Y\cap B_R.
\]
For a measure $\sigma$, we write
\[
    \sigma_R:=\sigma|_{B_R}
\]
for its restriction to $B_R$. For $p\in\R^d$, we denote by $\tau_p\nu$ the translation
of a measure $\nu$ by $p$, defined by
\[
    (\tau_p\nu)(H):=\nu(H-p)
\]
for Borel sets $H\subseteq\R^d$. We denote by $C_c(\R^d)$ the space of compactly supported continuous functions
and by $\schw$ the Schwartz class. Our Fourier transform convention is
\begin{equation}\label{eq:fourier-convention}
    \widehat f(\lambda)
    :=\int_{\R^d}f(x)e^{-2\pi i\langle x,\lambda\rangle}\,dx,
    \qquad f\in\schw.
\end{equation}
The space of tempered distributions is denoted by $\mathcal{S}'(\R^d)$. For
$T\in\mathcal{S}'(\R^d)$, its Fourier transform is defined by duality,
\[
    \widehat T(f):=T(\widehat f),
    \qquad f\in\schw.
\]
We refer to \cite[Chs.~8--9]{BaakeGrimmBook} for the basic facts on Fourier
analysis of measures and mathematical diffraction used below.

We write $A\lesssim B$ if $A\leq cB$ for a finite constant $c>0$ that does not
depend on the varying parameters involved (such as the radius $R$ or the points
$x$, $q$); a subscript, as in $\lesssim_\rho$, indicates that the constant may
depend on that parameter. In almost-sure statements, the constant may depend on
the realization unless otherwise stated. We write $A_R=O(B_R)$ if $|A_R|\leq C|B_R|$ for some
constant $C<\infty$ and all sufficiently large $R$, and
$A_R=o(B_R)$ if $A_R/B_R\to0$ as $R\to\infty$,
whenever $B_R\neq0$ for all sufficiently large $R$.
The same conventions on parameter dependence and almost-sure
statements apply.

\subsection{Radon measures, translation boundedness and convergence}

A complex Radon measure on $\R^d$ is a linear functional
$\mu:C_c(\R^d)\to\C$ such that, for every compact set $K\subseteq\R^d$, there
is a constant $C_K$ satisfying
\[
    |\mu(f)|\leq C_K\|f\|_\infty
\]
for every $f\in C_c(\R^d)$ supported in $K$. Its total variation is denoted by
$|\mu|$. In particular, for a purely atomic measure
\[
    \mu=\sum_{p\in X}\mu(p)\,\delta_p
\]
one has
\[
    |\mu|
    =\sum_{p\in X}|\mu(p)|\,\delta_p.
\]

A Radon measure $\mu$ is called \emph{translation bounded} (TB) if
\begin{equation}\label{eq:TB-def}
    \|\mu\|_{\mathrm{TB}}
    :=\sup_{x\in\R^d}|\mu|(B_1(x))<\infty.
\end{equation}
Equivalently, one may replace $B_1$ by any fixed compact neighborhood of the
origin; see \cite[Sec.~8.5]{BaakeGrimmBook}. A standard covering argument gives
a dimensional constant $C_d>0$ such that
\begin{equation}\label{eq:TB-poly-growth}
    |\mu|(B_R(x))
    \leq C_d\|\mu\|_{\mathrm{TB}}(1+R)^d,
    \qquad x\in\R^d,\quad R\geq0.
\end{equation}
In particular, every translation bounded measure has polynomial growth and
therefore defines a tempered distribution. For a purely atomic measure
$\mu=\sum_{p\in X}\mu(p)\,\delta_p$ we also write
$\|\mu\|_\infty:=\sup_{p\in X}|\mu(p)|$; note that $\|\mu\|_\infty\leq\|\mu\|_{\mathrm{TB}}$.

\medskip

We say that a family of Radon measures $\{\nu_R\}_{R\geq1}$ is
\emph{uniformly translation bounded} if
\[
    \sup_{R\geq1}\|\nu_R\|_{\mathrm{TB}}<\infty.
\]
By \eqref{eq:TB-poly-growth}, such a family satisfies a polynomial growth bound
uniformly in $R$. We say that $\nu_R$ converges \textit{vaguely} to $\nu$ if
\[
    \nu_R(f)\longrightarrow\nu(f),
    \qquad f\in C_c(\R^d).
\]
When the measures involved define tempered distributions, we say that
$\nu_R\to\nu$ \emph{in the distributional sense} if
\[
    \nu_R(f)\longrightarrow\nu(f),
    \qquad f\in\schw.
\]

\medskip

We will use the following elementary result, proved in \cite[Lemma~2.1]{nosotros}. 

\begin{lem}\label{lem:vague-to-distribution}
Let $\{\nu_R\}_{R\geq1}$ be a family of Radon measures for which there exist
constants $C>0$ and $n\in\N$ such that
\[
    |\nu_R|(B_k)\leq Ck^n,
    \qquad R\geq1,\quad k\geq1.
\]
Assume that $\nu_R\to\nu$ vaguely and that $\nu$ has polynomial growth. Then
$\nu_R\to\nu$ in $\mathcal{S}'(\R^d)$.
\end{lem}

In particular, if $\{\nu_R\}$ is uniformly translation bounded and converges
vaguely to a translation bounded measure $\nu$, then the convergence is automatically in the distributional sense. We will also use the following standard shell estimate for
translation bounded measures, see also
\cite[Lemma~2.2]{nosotros}.

\begin{lem}\label{lem:TB-shell}
There is a constant $C=C(d)>0$ such that, for every translation bounded
measure $\nu$, every $q\in\R^d$, and $R_2>R_1\geq0$,
\[
    |\nu|\bigl(B_{R_2}(q)\setminus B_{R_1}(q)\bigr)
    \leq
    C\|\nu\|_{\mathrm{TB}}
    (1+R_2)^{d-1}\max\{1,R_2-R_1\}.
\]
\end{lem}

\subsection{Convolution, reflection and Fourier transforms of measures}

If $\nu_1$ and $\nu_2$ are finite measures, their convolution is defined by
\[
    (\nu_1*\nu_2)(f)
    :=
    \int_{\R^d\times\R^d}
        f(x+y)\,d\nu_1(x)\,d\nu_2(y),
    \qquad f\in C_c(\R^d).
\]
More generally, if one of the measures is finite and the other is
translation bounded, their convolution is well defined and is again
translation bounded; see \cite[Ch.~8]{BaakeGrimmBook}.
For a function $f$ and a measure $\nu$, we use the reflected-conjugate
notation
\[
    \widetilde f(x):=\overline{f(-x)},
    \qquad
    \widetilde\nu(f):=\overline{\nu(\widetilde f)}.
\]
Thus,
\[
    \nu=\sum_{p}a_p\delta_p\quad \Longrightarrow \quad
    \widetilde\nu
    =\sum_p\overline{a_p}\,\delta_{-p}.
\]

For a finite measure $\sigma$, its Fourier--Stieltjes transform is the bounded
and uniformly continuous function
\begin{equation}\label{eq:finite-measure-FT}
    \widehat\sigma(\lambda)
    :=
    \int_{\R^d}
       e^{-2\pi i\langle x,\lambda\rangle}\,d\sigma(x).
\end{equation}
When $\sigma$ is regarded as a tempered distribution, its Fourier transform can be identified as an absolutely continuous measure having $\widehat\sigma$ as density with
respect to $m_d$. With our convention, $\widehat{\delta_0}=m_d$, and $\widehat{m_d}=\delta_0$
as identities of tempered distributions; see
\cite[Secs.~8.4--8.6]{BaakeGrimmBook}.
For a finite measure $\sigma$, the Fourier transform of
$\sigma*\widetilde\sigma$ is
\[
    \widehat{\sigma*\widetilde\sigma}(\lambda)
    =|\widehat\sigma(\lambda)|^2.
\]
Consequently, if $\sigma$ is a probability measure, then, as an identity of
tempered distributions,
\begin{equation}\label{eq:delta-minus-prob-FT}
    \widehat{\delta_0-\sigma*\widetilde\sigma}
    =
    \bigl(1-|\widehat\sigma|^2\bigr)m_d.
\end{equation}

We will also use the convolution theorem in the following form. If $\sigma$
is finite, $\nu$ is translation bounded, and $\widehat\nu$ is a measure, then
$\sigma*\nu$ is translation bounded and
\begin{equation}\label{eq:FT-conv-mult}
    \widehat{\sigma*\nu}
    =
    \widehat\sigma\,\widehat\nu,
\end{equation}
where the right-hand side denotes multiplication of the measure
$\widehat\nu$ by the bounded continuous function $\widehat\sigma$; see
\cite[Ch.~8]{BaakeGrimmBook}. More generally, for a measure $\nu$ and a bounded continuous function $h$, we write $h\nu$ for the measure defined by
\[
(h\nu)(f):=\int_{\mathbb R^d}f(x)h(x)\,d\nu(x).
\]

For an $\R^d$-valued random vector $\xi$, we denote its law by $\pe_\xi$ and
its characteristic function, with the Fourier convention
\eqref{eq:fourier-convention}, by
\begin{equation}\label{eq:characteristic-function}
    \varphi_\xi(\lambda)
    :=
    \widehat{\pe_\xi}(\lambda)
    =
    \e\bigl(e^{-2\pi i\langle\lambda,\xi\rangle}\bigr).
\end{equation}
Since $\widetilde{\pe_\xi}$ is the law of $-\xi$, one has
\begin{equation}\label{eq:diff-law-char}
    \widehat{\pe_\xi*\widetilde{\pe_\xi}}(\lambda)
    =
    |\varphi_\xi(\lambda)|^2.
\end{equation}

\subsection{Autocorrelation and diffraction}

For a Radon measure $\nu$ and $R>0$, its normalized finite-volume
autocorrelation is
\begin{equation}\label{eq:finite-auto-prelim}
    \gamma_\nu^{(R)}
    :=
    \frac{\nu_R*\widetilde{\nu_R}}{m_d(B_R)}.
\end{equation}

When the vague limit
\begin{equation}\label{eq:auto-def}
    \gamma_\nu
    :=
    \lim_{R\to\infty}\gamma_\nu^{(R)}
\end{equation}
exists, we call it the \emph{autocorrelation} of $\nu$. Since throughout this
paper we average over Euclidean balls, this is what is called the
\emph{natural autocorrelation} in \cite[Sec.~9.1.1]{BaakeGrimmBook}.

In particular, for our deterministic atomic measure $\mu$,
we have
\[
    \gamma_\mu^{(R)}
    =\frac{1}{m_d(B_R)}
    \sum_{p,q\in X_R}
    \mu(p)\overline{\mu(q)}\,\delta_{p-q}.
\]
For a translation bounded measure $\nu$, the family
$\{\gamma_\nu^{(R)}\}_{R\geq1}$ is uniformly translation bounded and vaguely
precompact. In particular, it always has vague accumulation points, and every
such accumulation point is itself translation bounded and positive definite;
see \cite[Prop.~9.1]{BaakeGrimmBook}.

This deterministic fact is worth keeping in mind in the random setting of the
present paper: the perturbed measure $\mu_\xi$ need not itself be translation
bounded, whereas we will prove that its normalized finite-volume
autocorrelations are nevertheless uniformly translation bounded almost surely.

Assume that the autocorrelation $\gamma_\nu$ defines a tempered distribution. This distribution is positive definite, and the Bochner--Schwartz theorem implies that its distributional Fourier transform is a positive tempered Radon measure. This measure, $\widehat{\gamma_\nu}$, is called the \emph{diffraction measure} of $\nu$; see \cite[Sec.~9.1.2]{BaakeGrimmBook} and \cite{Hof}.
Its Lebesgue decomposition is written as
\[
    \widehat{\gamma_\nu}
    =
    \bigl(\widehat{\gamma_\nu}\bigr)_{\mathrm{pp}}
    +
    \bigl(\widehat{\gamma_\nu}\bigr)_{\mathrm{sc}}
    +
    \bigl(\widehat{\gamma_\nu}\bigr)_{\mathrm{ac}},
\]
where the three terms are respectively the pure point, singular continuous,
and absolutely continuous components, with singularity and absolute
continuity understood with respect to Lebesgue measure on $\mathbb R^d$. The pure point part describes the
Bragg peaks, while the continuous components correspond to diffuse
scattering; see \cite[Sec.~9.1.2]{BaakeGrimmBook}.

The existence of the diffraction measure must be distinguished from its computation as a limit of finite-volume approximants. Indeed, once the random autocorrelation $\gamma_{\mu_\xi}$
exists, its Fourier transform is well defined and gives the
diffraction measure $\widehat{\gamma_{\mu_\xi}}$. However, since Fourier transforms are understood in the sense of tempered
distributions, as in \cite{BaakeGrimmBook,Hof}, vague convergence
alone does not justify passing to the Fourier transforms. Thus, to obtain, in addition, 
\[
    \lim_{R\to\infty}\gamma_{\mu_\xi}^{(R)}(\widehat f)
    =\widehat{\gamma_{\mu_\xi}}(f),
    \qquad f\in\schw,
\]
we strengthen the vague convergence of the finite-volume autocorrelations to convergence in \(\mathcal S'(\mathbb R^d)\). Theorem~\ref{thm:uniform-TB}, together with Lemma~\ref{lem:vague-to-distribution}, provides this upgrade, and continuity of the Fourier transform then yields the identity.

\subsection{Random measures and probabilistic notation}
\label{subsec:probabilistic-preliminaries}

By a random finite measure $\eta$ we mean a
measure-valued random element such that $\eta(H)$
is measurable for every Borel set $H\subseteq\R^d$. For a random finite complex measure $\eta$ satisfying
$\mathbb E|\eta|(\R^d)<\infty$, its mean measure
$\mathbb E\eta$ is defined by
\[
    (\mathbb E\eta)(H):=\mathbb E[\eta(H)],
\]
for every Borel set $H\subseteq\R^d$. In particular, $ |\mathbb E\eta|(\R^d)\leq \mathbb E|\eta|(\R^d)$, and 
for every bounded Borel function $f$,
\[
    (\mathbb E\eta)(f)=\mathbb E[\eta(f)].
\]

If $\eta$ is positive, the identity for test functions also
holds for nonnegative Borel functions, by Tonelli's theorem.

For a random element $Y$, we denote by $\mathcal F_Y$
the sigma-algebra generated by $Y$, and by
$\mathcal F_{(Y_1,\ldots,Y_k)}$ the sigma-algebra generated
by the corresponding tuple. We write
$\mathbb E[Z\mid\mathcal F_Y]$ for the conditional
expectation of $Z$ given $\mathcal F_Y$.

\section{Random perturbations of TB purely atomic measures}\label{randomperturbTB}

The purpose of this section is to establish the local properties of the randomly perturbed measure that will be needed from the next section onward in the study of its autocorrelation. Recall that 
\[
 \mu = \sum_{p\in X} \mu(p) \delta_p
\]
is a translation bounded (TB) measure with locally finite support $X\subseteq \R^d$, and that $\{\xi_p\}_{p\in X}$ is an i.i.d.\ family of copies of the random vector $\xi$. The existence of $\gamma_\mu$ is not necessary in this section. Consider 
\[
\mu_\xi=\sum_{p\in X}\mu(p)\delta_{p+\xi_p} \,.
\]
We begin by showing that $\mu_\xi$ is almost surely a Radon measure and it is translation bounded in mean. Since different atoms may be displaced to the same location, cancellations may occur, so the total variation of $\mu_\xi$ need not equal the sum of the absolute values of the displaced coefficients. It is therefore convenient to introduce the positive Borel measure
\[
 \rho_\xi:=\sum_{p\in X}|\mu(p)|\,\delta_{p+\xi_p} \,.
\]
These properties are established below in a more general
setting that also covers random marks, so that the same
result can be used in Sections~\ref{averaging}
and~\ref{generalcase}.

\begin{lem}\label{lem:local-mark-mass}
Let $(\alpha_p)_{p\in X}$ be identically distributed random
finite positive measures of mass at most $C$ almost surely,
and set $A_p=\tau_p\alpha_p$, where \(\tau_p\) is the translation by \(p\). Then the positive measure
\[
    \rho_A:=\sum_{p\in X}|\mu(p)|A_p
\]
is almost surely a Radon measure. Moreover, for every
bounded Borel set $H\subseteq\R^d$,
\[
    \mathbb E\rho_A(H)
    =\int|\mu|(H-u)\,d\bar\alpha(u)
    \le C\sup_{u\in\R^d}|\mu|(H-u),
\]
where $\bar\alpha=\mathbb E\alpha_p$.

If $(a_p)_{p\in X}$ is a family of random finite complex
measures satisfying $|a_p|\le A_p$ almost surely, then
\[
    \nu:=\sum_{p\in X}\mu(p)a_p
\]
is almost surely a Radon measure and $|\nu|\le\rho_A$.
\end{lem}

\begin{proof}
Since $A_p(H)=\alpha_p(H-p)$, Tonelli's theorem gives
\[
    \mathbb E\rho_A(H)
    =\int\sum_{p\in X}|\mu(p)|
          \mathbbm1_H(p+u)\,d\bar\alpha(u)
    =\int|\mu|(H-u)\,d\bar\alpha(u).
\]
The asserted bound follows from
$\bar\alpha(\R^d)\le C$ and the translation boundedness
of $\mu$. Taking $H=B_n$ and intersecting over the positive integers
gives a single event of probability one on which $\rho_A$
is locally finite. On this event, the series defining
$\nu$ converges absolutely in total variation on compact
sets, and $|\nu|\le\rho_A$.
\end{proof}

\begin{obs*}
Independence across sites is not needed in
Lemma~\ref{lem:local-mark-mass}.
\end{obs*}

\begin{prop}\label{radon}
    Almost surely, both $\rho_\xi$ and $\mu_\xi$ are Radon measures and satisfy $|\mu_\xi| \leq \rho_\xi$.
\end{prop}

\begin{proof}
Apply Lemma~\ref{lem:local-mark-mass} with
$A_p=a_p=\delta_{p+\xi_p}$ and $C=1$.
\end{proof}

The mean bound in Lemma~\ref{lem:local-mark-mass}
also gives the following corollary.

\begin{cor}\label{tbinmean}
   For each compact $K$ \[\sup_{x\in \mathbb{R}^d} \e(\rho_\xi(K+x)) < \infty \,.\] 
\end{cor}
Taking $K = B_1$ in this corollary, we get that $\mu_\xi$ is \textit{TB in mean}. Despite this, it is not difficult to see that, in general, $\mu_\xi$ fails to be TB. Indeed, we have the following example:

\begin{exa}\label{ex:lattice-gaussian-not-TB}
Let \(\mu=\delta_{\Z}\), and suppose the displacement law charges
every open interval. Fix \(n\ge1\). For each \(k\ge0\), the \(n\)
atoms at \(kn,\dots,kn+n-1\) all land in \((kn-1,kn+1)\) with the
same positive probability. These events are independent because
the blocks are disjoint. Almost surely infinitely many occur.
Intersecting over \(n\) shows that \(\mu_\xi\) is almost surely
not translation bounded.
\end{exa}

\subsection{Local mass after displacement}

For a translation bounded measure, the total variation
over translates of a fixed cube is uniformly bounded.
As the preceding example shows, this property may fail
for $\mu_\xi$. This motivates the following weaker notion
of translation boundedness.

Write \(Q=[-1/2,1/2)^d\). A Radon measure \(\nu\) is
\emph{quasi-translation bounded} (quasi-TB) if for every \(\delta>0\),
\begin{equation}\label{quasi}
 \lim_{\|q\|\to\infty}\frac{|\nu|(Q+q)}{\|q\|^\delta}=0 .
\end{equation}
It suffices to take \(q\in\Z^d\), since any translate of \(Q\) is
covered by a bounded number of lattice translates with nearby
centers. For \(R>e\), write
\[
 \ell(R)=\frac{\log R}{\log\log R}.
\]

\begin{teo}\label{emeR}
Let \(D>d\,\|\mu\|_\infty\). Almost surely,
\(\rho_\xi(Q+q)\le D\,\ell(\|q\|)\) for all
\(q\in\Z^d\) of sufficiently large norm. In particular, \(\rho_\xi\)
and \(\mu_\xi\) are quasi-TB.
\end{teo}

\begin{obs*}
Proposition~\ref{radon} and Theorem~\ref{emeR} apply to arbitrary
displacement laws, including those with infinite first moment.
\end{obs*}

The local mass estimate follows from a bound for sums of
bounded independent variables that also applies to the marked
measures in Section~\ref{generalcase}.

\begin{lem}\label{lem:lattice-growth}
Let \(J\) be countable. For each \(q\in\Z^d\), suppose
\((W_{p,q})_{p\in J}\) are independent random variables with
\(0\le W_{p,q}\le b\) for a common constant \(b>0\).
If \(\sup_{q\in\Z^d}\sum_{p\in J}\mathbb E W_{p,q}<\infty\), then
for every \(a>d\), almost surely,
\[
 \sum_{p\in J} W_{p,q}\le ab\,\ell(\|q\|)
 \qquad\text{for all \(q\in\Z^d\) of sufficiently large norm}.
\]
\end{lem}

\begin{proof}
Write \(T_q=\sum_{p\in J}W_{p,q}\) and fix \(m>0\) with
\(\mathbb E T_q\le m\) for every \(q\). For a real parameter \(r>1\),
convexity gives \(r^u\le1+(r-1)u\) for \(0\le u\le1\). Thus for every
finite \(F\subset J\), independence and \(1+x\le e^x\) give
\[
 \mathbb E \left(r^{\sum_{p\in F} W_{p,q}/b}\right)
 \le \prod_{p\in F}\left(1+\frac{r-1}{b}\mathbb E W_{p,q}\right)
 \le \exp\!\left(\frac{m(r-1)}{b}\right).
\]
Monotone convergence gives the same bound for \(\mathbb E r^{T_q/b}\).
For \(t>m\), Markov's inequality with \(r=t/m\) yields
\[
 \mathbb P(T_q\ge t)
 \le \exp\!\left(\frac{t-m-t\log(t/m)}{b}\right)
 \le \left(\frac{em}{t}\right)^{t/b}.
\]
Given \(a>d\), take \(t=ba\,\ell(\|q\|)\) at sufficiently large
\(\|q\|\). The last bound is then \(\|q\|^{-a+o(1)}\), which is
summable over \(\Z^d\). Borel--Cantelli gives the conclusion.
\end{proof}

\begin{proof}[Proof of Theorem~\ref{emeR}]
The variables \(|\mu(p)|\mathbbm1_{Q+q}(p+\xi_p)\) are independent
over \(p\) and bounded by \(\|\mu\|_\infty\). Their sum is
\(\rho_\xi(Q+q)\). Its expectation is bounded uniformly in \(q\)
by Corollary~\ref{tbinmean} applied to the closure of \(Q\).
Lemma~\ref{lem:lattice-growth}, applied with \(b=\|\mu\|_\infty\)
and \(a=D/b>d\), gives the asserted bound. The
quasi-TB conclusion for \(\mu_\xi\) follows from
\(|\mu_\xi|\le\rho_\xi\).
\end{proof}

The following result shows that the order $\ell(\|q\|)$ in Theorem~\ref{emeR} is attained, even in the simplest case of a lattice whose points are perturbed by using a displacement law with finite first moment.

\begin{prop}\label{ex:sharp-local-growth}
For every $d\geq1$, there exists a non-atomic displacement law with finite first moment such that, for $\mu=\delta_{\Z^d}$, almost surely
\[
    0<\limsup_{\substack{q\in\Z^d\\ \|q\|\to\infty}}
    \frac{\mu_\xi(Q+q)}{\ell(\|q\|)}<\infty.
\]
In particular, the order of growth in Theorem~\ref{emeR} cannot be improved uniformly over the displacement laws.
\end{prop}

\begin{proof}
Fix $\beta>d+1$, and let $\pe_\xi$ be the probability measure on $\R^d$ whose density is
\[
    \psi_\beta(x)=c_\beta\,(1+\|x\|)^{-\beta},
\]
where $c_\beta$ is the normalizing constant. This law is non-atomic, and it has finite first moment because $\beta-d>1$. Let $e_1$ denote the first vector of the canonical basis of $\R^d$. Given two integers $k\geq1$ and $0\leq j<k$, every $x\in Q-je_1$ satisfies $\|x\|\leq k+\sqrt d$. Since $m_d(Q)=1$, this implies that
\begin{equation}\label{eq:sharp-one-site}
    \pe(\xi\in Q-je_1)\geq b\,k^{-\beta},
    \peso{where} b:=c_\beta\,(2+\sqrt d)^{-\beta}.
\end{equation}
Now, fix $a\in(0,1/\beta)$, let $k_m:=\lfloor a\,\ell(m)\rfloor$, and take $m_0\geq16$ such that $k_m\geq1$ for every $m\geq m_0$. Starting with $s_{m_0}:=1$, define recursively $s_{m+1}:=s_m+k_m$, and consider the pairwise disjoint blocks of lattice points $\cI_m:=\{(s_m+j)\,e_1:\ 0\leq j<k_m\}$. Let $F_m$ be the event that all the atoms located at the points of $\cI_m$ are displaced into the cube $Q+s_me_1$. Since the atom located at $(s_m+j)e_1$ lands in this cube if and only if its displacement belongs to $Q-je_1$, by the independence of the displacements and~\eqref{eq:sharp-one-site} we get that
\[
    \pe(F_m)\geq\bigl(b\,k_m^{-\beta}\bigr)^{k_m}
    =\exp\bigl(-\beta\,k_m\log k_m-k_m\log(1/b)\bigr)
    \geq m^{-a\beta-o(1)},
\]
where in the last inequality we have used that $k_m\log k_m\leq a\log m$ and $k_m=o(\log m)$. As $a\beta<1$, the series $\sum_m\pe(F_m)$ diverges. Since the events $F_m$ are independent, because they involve disjoint blocks of sites, by the second Borel--Cantelli lemma almost surely infinitely many of them occur.

On the event $F_m$, we have that $\mu_\xi(Q+s_me_1)\geq k_m$. On the other hand, since $1\leq k_i\leq\log i$ for every $i\geq m_0$, we have that $m-m_0\leq s_m\leq1+m\log m$. Hence, $\ell(s_m)/\ell(m)\to1$, and therefore $k_m/\ell(s_m)\to a$. In consequence, almost surely the upper limit of the statement is at least $a>0$. Finally, it is finite by Theorem~\ref{emeR}, because in this case $\rho_\xi=\mu_\xi$.
\end{proof}

\subsection{Boundary-crossing estimates}

We now estimate the total weight of the points
that enter or leave $B_R$ under the perturbation.
Let \(X_R=X\cap B_R\) and
\(X_{R,\xi}=\{p\in X:p+\xi_p\in B_R\}\). The total weight of the
sites that enter or leave the window is
\[
\begin{aligned}
 D_R^\xi
 &=|\mu|(X_R\mathbin{\triangle}X_{R,\xi})=\sum_{p\in X}|\mu(p)|
    \int\bigl|\mathbbm1_{B_R}(x)-\mathbbm1_{B_R}(p)\bigr|
       \,d\delta_{p+\xi_p}(x).
\end{aligned}
\]
These estimates will be a main ingredient in the next
section where, among other things, we show that the random truncation
indexed by $X_{R,\xi}$ can be asymptotically replaced
by the deterministic truncation indexed by $X_R$,
without any moment assumption on the displacements.
That window comparison improves the corresponding estimates in \cite[Lemma~2]{Yakir2022}, \cite[Lemma~4.1]{PetracheViera2023}
and \cite[Lemma~4.2]{nosotros}.

\medskip

\begin{prop}\label{entranysalen}
For any displacement law, \(D_R^\xi=o(R^d)\) almost surely.
\end{prop}

We will use the following summability criterion here and in the
pair estimates of the next section. 

\begin{lem}\label{lem:quadratic-summability}
Let \((T_n)\) be square-integrable random variables on a common
probability space and \(a_n>0\). If
\(\sum_n\operatorname{Var}(T_n)/a_n^2<\infty\), then
\[
 \sum_n\left|\frac{T_n-\mathbb E T_n}{a_n}\right|^2<\infty
 \qquad\text{almost surely}.
\]
If also \(\mathbb E T_n=o(a_n)\), then \(T_n/a_n\to0\) almost
surely.
\end{lem}

\begin{proof}
By Tonelli, the series of squares has expectation
\(\sum_n\operatorname{Var}(T_n)/a_n^2\), so it is finite almost
surely. Its terms tend to zero, and adding \(\mathbb E T_n/a_n\)
gives the last assertion.
\end{proof}

To treat both displacements and marks at the same time, let \(A_p=\tau_p\alpha_p\),
where the \(\alpha_p\) are i.i.d.\ positive measures of mass at
most \(C\) almost surely. For the family \(A=(A_p)_{p\in X}\), the crossing mass is
\[
 D_R^A=\sum_{p\in X}|\mu(p)|
 \int\bigl|\mathbbm1_{B_R}(x)-\mathbbm1_{B_R}(p)\bigr|\,dA_p(x).
\]
Equivalently, separating the mass leaving the window
from that entering it, we have
\[
 D_R^A=
 \sum_{p\in X_R}|\mu(p)|A_p(B_R^c)
 +\sum_{p\in X\setminus X_R}|\mu(p)|A_p(B_R).
\]
Taking \(A_p=\delta_{p+\xi_p}\) recovers the integral formula for \(D_R^\xi\).
Thus Proposition~\ref{entranysalen} is a special case of the following result.

\begin{lem}\label{coronas}
The crossing mass satisfies \(\mathbb E D_R^A=o(R^d)\) and
\(D_R^A=o(R^d)\) almost surely.
\end{lem}
\begin{proof}
Let \(\bar\alpha=\mathbb E\alpha_p\). Setting \(x=p+u\) in the
integral and using Tonelli gives
\[
 \mathbb E D_R^A
 =\int|\mu|\bigl(B_R\mathbin{\triangle}(B_R-u)\bigr)
                  \,d\bar\alpha(u).
\]
Each ball in the symmetric difference has mass \(O(R^d)\)
uniformly in its center. For \(0<\|u\|<R\), the difference is
also contained in the union of two inner shells of width
\(\|u\|\), one in each ball. By Lemma~\ref{lem:TB-shell},
\[
 |\mu|\bigl(B_R\mathbin{\triangle}(B_R-u)\bigr)
 \lesssim R^{d-1}\min\{R,1+\|u\|\},\qquad R\ge1.
\]
Consequently the integrand divided by \(R^d\) tends to zero
for each fixed \(u\) and is bounded uniformly in \(u\).
Dominated convergence gives the mean estimate.

The site contributions to \(D_R^A\) are independent and bounded
by \(C\|\mu\|_\infty\). Bounding the variance of each contribution
by this constant times its mean, and passing from finite sums by
monotone convergence, gives
\(\operatorname{Var}(D_R^A)\le C\|\mu\|_\infty\,\mathbb E D_R^A\).
At \(R=n^2\) the variance bound gives
\(\operatorname{Var}(D_{n^2}^A)=O(n^{2d})\). Since
\(\sum_n n^{-2d}<\infty\), Lemma~\ref{lem:quadratic-summability}
gives \(D_{n^2}^A=o(n^{2d})\) almost surely.

If \(n^2\le R\le(n+1)^2\), a crossing whose original site is
outside the intervening shell is counted at one of the two
endpoint radii. Hence
\[
 D_R^A\le D_{n^2}^A+D_{(n+1)^2}^A
              +C|\mu|\bigl(B_{(n+1)^2}\setminus B_{n^2}\bigr).
\]
The shell has mass \(O(n^{2d-1})\) by Lemma~\ref{lem:TB-shell},
and \((n+1)^{2d}\sim n^{2d}\). Thus the right-hand side is
\(o(n^{2d})\), which gives the limit along all radii.
\end{proof}

\begin{cor}\label{mejora}
Almost surely, \(\rho_\xi(B_R)\lesssim R^d\) as \(R\to\infty\).
\end{cor}

\begin{proof}
Every displaced atom in \(B_R\) either comes from a site in \(B_R\)
or has crossed the boundary, so
\(\rho_\xi(B_R)\le|\mu|(B_R)+D_R^\xi\). The claim follows from
\eqref{eq:TB-poly-growth} and Proposition~\ref{entranysalen}.
\end{proof}

\section{Window comparison and centering}\label{averaging}

A marked configuration gives two ways to form a finite-volume
convolution. The observed window restricts the mass at the final
positions to \(B_R\). The original window selects the sites in
\(X_R=X\cap B_R\) and retains the full mark at each of these sites.
We compare the two convolutions at the volume scale, then estimate
the fluctuations of the original-window convolution about its mean.

Let \((\alpha_p,\beta_p)_{p\in X}\) be i.i.d.\ pairs of positive
measures, each of mass at most \(C\) almost surely. The two marks
at a site may be dependent. Write \(\Phi_p=\tau_p\alpha_p\) and
\(\Psi_p=\tau_p\beta_p\) for their translates to the final positions.

Let \(a_p,b_p\) be random finite complex measures on final positions,
with \(|a_p|\le \Phi_p\) and \(|b_p|\le \Psi_p\) almost surely. Restricting
the final positions to \(B_R\) gives the measures
\[
 \sigma_R=\sum_{p\in X}\mu(p)\,a_p|_{B_R},
 \qquad
 \lambda_R=\sum_{p\in X}\mu(p)\,b_p|_{B_R}.
\]
A prime indicates that the window selects the original sites. Each
selected mark is retained in full, including any mass outside \(B_R\).
\[
 \sigma'_R=\sum_{p\in X_R}\mu(p)a_p,
 \qquad
 \lambda'_R=\sum_{p\in X_R}\mu(p)b_p.
\]

Lemma~\ref{lem:local-mark-mass} gives almost-sure convergence
in total variation of the observed-window sums
$\sigma_R$ and $\lambda_R$. The primed sums are finite
because $X_R$ is finite.

\begin{teo}\label{soportes}
On one event of probability one,
\[
 \frac{\sigma_R*\widetilde{\lambda_R}
         -\sigma'_R*\widetilde{\lambda'_R}}
      {m_d(B_R)}
 \xrightarrow[R\to\infty]{\mathrm{vague}}0 .
\]
\end{teo}
In the original window, the convolution is indexed by pairs of
original sites. For \(f\in C_c(\R^d)\), its evaluation is the random
sum
\[
 (\sigma'_R*\widetilde{\lambda'_R})(f)
 =\sum_{p,q\in X_R}\mu(p)\overline{\mu(q)}\,Z_{p,q}(f),
 \qquad
 Z_{p,q}(f)=\iint f(x-y)\,da_p(x)\,d\overline{b_q}(y).
\]
Here $\overline{b_q}$ denotes the conjugate measure of
$b_q$, defined by $\overline{b_q}(H)=\overline{b_q(H)}$
for Borel sets $H\subseteq\R^d$.

The next theorem is a strong law of large numbers for these pair
sums. As in Kolmogorov's strong law for independent, not identically
distributed variables, it is stated in centered form. It asserts that
the fluctuations are negligible at the volume scale, but it does not
require the normalized means to converge. Whether they converge is a
deterministic question, which we address in
Section~\ref{randomautocorrelation}.

\begin{teo}\label{lem:parent-centering}
Suppose the quadruples \((a_p,b_p,\Phi_p,\Psi_p)\) are independent
across sites. Then on one event of probability one, for every
\(f\in C_c(\R^d)\),
\[
 \frac{(\sigma'_R*\widetilde{\lambda'_R})(f)
       -\mathbb E[(\sigma'_R*
                    \widetilde{\lambda'_R})(f)]}
      {m_d(B_R)}
 \xrightarrow[R\to\infty]{}0 .
\]
\end{teo}

The proofs of both theorems rely on two auxiliary estimates.
The first controls the mass that can pair with a given mark
within a fixed distance, after conditioning on the data at
its original site. Together with the boundary-crossing
estimate of Lemma~\ref{coronas}, this allows us to control
the error in the window comparison. The second is a variance
bound for pair sums, which exploits the independence of
terms indexed by disjoint pairs of sites. We first state
these estimates, then prove the two theorems, and finally
give the proofs of the auxiliary lemmas.

\subsection{Local mass and variance bounds}

Recall that each positive mark \(\Phi_p\) or \(\Psi_p\) has total mass
at most \(C\) almost surely, with the same bound for all sites.
For a radius \(\rho\ge1\), set
\[
 m_\rho=\sup_{x\in\R^d}|\mu|(B_\rho(x)),
\]
and recall that \(\|\mu\|_\infty=\sup_{p\in X}|\mu(p)|\).
The mass of pairs with final positions within distance \(\rho\) is
given by the positive kernel
\begin{equation}\label{eq:proximity-kernel}
 V_{p,q}^{(\rho)}
 =\iint\mathbbm1_{B_\rho(x)}(y)\,d\Phi_p(x)\,d\Psi_q(y),
 \qquad 0\le V_{p,q}^{(\rho)}\le C^2.
\end{equation}
Its weighted row with first index \(p\) integrates, against \(\Phi_p\),
the nearby mass contributed by the \(\Psi_q\). Its column with second
index \(p\) fixes \(\Psi_p\).

\begin{lem}\label{lem:local-pair-bounds}
For every
\(\rho\ge1\),
\begin{equation}\label{eq:mean-local-mark-mass}
 \sum_q|\mu(q)|\,\mathbb E \Psi_q(B_\rho(x))\le C m_\rho,\qquad x\in\R^d,
\end{equation}
and the same bound holds with \(\Phi_q\) in place of \(\Psi_q\). For each
\(p\in X\), almost surely,
\begin{equation}\label{eq:proximity-rows}
\begin{aligned}
 \mathbb E\!\left[\sum_q|\mu(q)|V_{p,q}^{(\rho)}
                 \,\middle|\,\mathcal{F}_{(\Phi_p,\Psi_p)}\right]
 &\le C^2(m_\rho+\|\mu\|_\infty),\\
 \mathbb E\!\left[\sum_q|\mu(q)|V_{q,p}^{(\rho)}
                 \,\middle|\,\mathcal{F}_{(\Phi_p,\Psi_p)}\right]
 &\le C^2(m_\rho+\|\mu\|_\infty).
\end{aligned}
\end{equation}
Under the independence hypothesis of Theorem~\ref{lem:parent-centering},
the same bounds hold when conditioning on \(\mathcal{F}_{(a_p,b_p,\Phi_p,\Psi_p)}\),
the sigma-algebra generated by the four measures at \(p\).
\end{lem}

\medskip

For the variance bound, let \((Y_p)_{p\in X}\) be independent
random elements. Given
\(\mathcal{F}_{(Y_p,Y_q)}\)-measurable kernels
\(Z_{p,q},G_{p,q}\) with \(|Z_{p,q}|\le G_{p,q}\le L\), and
deterministic coefficients \(c_p\) satisfying \(|c_p|\le|\mu(p)|\),
write
\[
 S_F=\sum_{p,q\in F}c_p\overline{c_q}Z_{p,q},
 \qquad
 S_F^{+}=\sum_{p,q\in F}|\mu(p)|\,|\mu(q)|G_{p,q}
 \quad(F\subset X\text{ finite}).
\]
Write \(S^{+}\) for the full positive sum over \(X\times X\).

\begin{lem}\label{lem:window-variance}
Suppose the conditional means of the weighted rows and columns
of \(G\), given \(\mathcal{F}_{Y_p}\) at the fixed site \(p\), are
bounded by \(d_0\). Then
\[
 \operatorname{Var}(S_F)\le 8\|\mu\|_\infty(\|\mu\|_\infty L+d_0)\,\mathbb E S_F^{+}
 \qquad(F\subset X\text{ finite}).
\]
If \(\mathbb E S^{+}<\infty\), then \(S^{+}\) is square-integrable and
\(\operatorname{Var}(S^{+})\le8\|\mu\|_\infty(\|\mu\|_\infty L+d_0)\,\mathbb E S^{+}\).
\end{lem}

\subsection{Proofs of the window comparison and centering}

\begin{proof}[Proof of Theorem~\ref{soportes}]
Fix \(f\in C_c(\R^d)\) supported in \(B_\rho\), with \(\rho\ge1\)
an integer. By Lemma~\ref{lem:local-mark-mass}, the observed-window
sums $\sigma_R$ and $\lambda_R$ converge absolutely
in total variation, since \(\sum_p|\mu(p)|\,|a_p|(B_R)\) and
\(\sum_p|\mu(p)|\,|b_p|(B_R)\) are bounded by the locally finite
dominating sums. The series for the observed convolution is also
absolutely convergent. Its absolute sum is bounded by the product
of these two sums times \(\|f\|_\infty\). The convolution in the
original window involves only finitely many sites.

Write \(\Delta_R(f)\) for the difference of the two convolutions
evaluated on \(f\). Expanding it gives
\[
\begin{aligned}
 \Delta_R(f)
 &=\sum_{p,q\in X}\mu(p)\overline{\mu(q)}\iint f(x-y)\times
 \bigl[\mathbbm1_{B_R}(x)\mathbbm1_{B_R}(y)
       -\mathbbm1_{B_R}(p)\mathbbm1_{B_R}(q)\bigr]\,
       da_p(x)\,d\overline{b_q}(y).
\end{aligned}
\]
In the first
product of indicators the window acts on the final positions
\(x,y\), and in the second it acts on the original sites \(p,q\).
The difference of the products in absolute value is bounded by
the sum of the two individual discrepancies, one for \((p,x)\)
and one for \((q,y)\).

Take \(r_n=n^2\) and \(s_n=r_{n+1}\). For \(r_n\le R\le s_n\),
the discrepancy \(|\mathbbm1_{B_R}(x)-\mathbbm1_{B_R}(p)|\) is bounded by
\[
\begin{aligned}
 h_n(p,x):={}&|\mathbbm1_{B_{r_n}}(x)-\mathbbm1_{B_{r_n}}(p)|
           +|\mathbbm1_{B_{s_n}}(x)-\mathbbm1_{B_{s_n}}(p)|+\mathbbm1_{B_{s_n}\setminus B_{r_n}}(p).
\end{aligned}
\]
Thus \(0\le h_n\le3\). Set
\[
 Z_{p,q}^{(n)}
 =\iint h_n(p,x)\,\mathbbm1_{B_\rho(x)}(y)\,d\Phi_p(x)\,d\Psi_q(y),
 \qquad
 S_n=\sum_{p,q\in X}|\mu(p)|\,|\mu(q)|Z_{p,q}^{(n)}.
\]
Let \(S'_n\) be the same sum with the two families exchanged.
Taking absolute values in the expansion of \(\Delta_R(f)\), use
the bounds \(|a_p|\le \Phi_p\), \(|\overline{b_q}|=|b_q|\le \Psi_q\),
and \(|f(x-y)|\le\|f\|_\infty\mathbbm1_{B_\rho(x)}(y)\).
The term with \(h_n(p,x)\) gives \(S_n\). The term with \(h_n(q,y)\)
gives \(S'_n\) upon interchanging \((p,x)\) and \((q,y)\).
Consequently,
\begin{equation}\label{eq:window-block-error}
 |\Delta_R(f)|\le\|f\|_\infty(S_n+S'_n),
 \qquad r_n\le R\le s_n.
\end{equation}

Since $Z_{p,q}^{(n)}\le3V_{p,q}^{(\rho)}$, where
$V_{p,q}^{(\rho)}$ is the proximity kernel defined
in \eqref{eq:proximity-kernel}, the conditional means of its weighted rows and columns are uniformly bounded by
\eqref{eq:proximity-rows}, applied with the pair \((\Phi_p,\Psi_p)\)
at the fixed site.
For the mean of \(S_n\) we retain \(h_n\) in the row estimate.
Given \(\mathcal{F}_{(\Phi_p,\Psi_p)}\), the contribution from \(q\ne p\)
uses \eqref{eq:mean-local-mark-mass}, integrated
against \(h_n(p,x)\,d\Phi_p(x)\). The diagonal uses
\(\Psi_p(B_\rho(x))\le C\) and \(|\mu(p)|\le \|\mu\|_\infty\). Together these give
\[
 \mathbb E\!\left[\sum_q|\mu(q)|Z_{p,q}^{(n)}
                   \,\middle|\,\mathcal{F}_{(\Phi_p,\Psi_p)}\right]
 \le C(m_\rho+\|\mu\|_\infty)\int h_n(p,x)\,d\Phi_p(x).
\]
Taking expectations and summing over \(p\) with weight
\(|\mu(p)|\), we get
\[
 \mathbb E S_n
 \le C(m_\rho+\|\mu\|_\infty)\bigl(\mathbb E D_{r_n}^\Phi+\mathbb E D_{s_n}^\Phi
                 +C|\mu|(B_{s_n}\setminus B_{r_n})\bigr)
 =o(r_n^d).
\]
Here \(D_R^\Phi\) is the crossing mass for the family
\((\Phi_p)_{p\in X}\). The two mean crossing terms are controlled by
Lemma~\ref{coronas} and the shell by Lemma~\ref{lem:TB-shell},
using \(s_n/r_n\to1\). The same estimate holds for \(S'_n\).

Taking \(Y_p=(\Phi_p,\Psi_p)\) in Lemma~\ref{lem:window-variance} gives
\(\operatorname{Var}(S_n)\lesssim_\rho\mathbb E S_n\) and the
corresponding bound for \(S'_n\).
By Lemma~\ref{lem:quadratic-summability}, both sums are
\(o(r_n^d)\) almost surely, since \(\sum_n r_n^{-d}<\infty\).
Dividing \eqref{eq:window-block-error} by \(m_d(B_R)\) gives
the limit throughout the interval \([r_n,s_n]\). For each \(\rho\)
the event depends only on the two positive sums, so it works for
every test function supported in \(B_\rho\). Intersect over integer \(\rho\).
\end{proof}

\begin{proof}[Proof of Theorem~\ref{lem:parent-centering}]
For \(f\) supported in \(B_\rho\), the kernel of the convolution is
\[
 Z_{p,q}(f)=\iint f(x-y)\,da_p(x)\,d\overline{b_q}(y),
 \qquad |Z_{p,q}(f)|\le\|f\|_\infty V_{p,q}^{(\rho)}.
\]
Fix \(\rho\) and take \(Y_p=(a_p,b_p,\Phi_p,\Psi_p)\)
in Lemma~\ref{lem:window-variance}.
First consider any \(\mathcal{F}_{(Y_p,Y_q)}\)-measurable
kernel \(Z_{p,q}\) bounded in absolute value by \(V_{p,q}^{(\rho)}\).
With \(|c_p|\le|\mu(p)|\), write
\[
 T_R=\sum_{p,q\in X_R}c_p\overline{c_q}(Z_{p,q}-\mathbb E Z_{p,q}),
 \qquad
 W_R=\sum_{p,q\in X_R}|\mu(p)|\,|\mu(q)|V_{p,q}^{(\rho)}.
\]
The conditional bounds of Lemma~\ref{lem:local-pair-bounds}
apply to these quadruples. The row bound gives
\(\mathbb E W_R\lesssim_\rho|\mu|(B_R)\), so
Lemma~\ref{lem:window-variance} gives
\(\mathbb E|T_R|^2=O(R^d)\).
With \(r_n=n^2\), Lemma~\ref{lem:quadratic-summability}
then gives \(T_{r_n}=o(r_n^d)\) almost surely.

Set \(s_n=r_{n+1}\) and \(H_n=W_{s_n}-W_{r_n}\). Each term of
\(H_n\) has at least one original site in
\(B_{s_n}\setminus B_{r_n}\). Summing the row and column bounds
over those sites gives
\[
 \mathbb E H_n\lesssim_\rho|\mu|(B_{s_n}\setminus B_{r_n})=o(r_n^d).
\]
Restricting the positive kernel to these pairs preserves the
conditional bounds in \eqref{eq:proximity-rows}. Hence
\(\operatorname{Var}(H_n)\lesssim_\rho\mathbb E H_n\), and the
summability lemma yields \(H_n=o(r_n^d)\) almost surely. For
\(r_n\le R\le s_n\), the centered increment is bounded by
\[
 |T_R-T_{r_n}|\le H_n+\mathbb E H_n.
\]
This proves \(T_R/m_d(B_R)\to0\) along all radii.

First apply the scalar conclusion with \(Z_{p,q}=V_{p,q}^{(\rho)}\)
and \(c_p=|\mu(p)|\). The centered sum in this case is
\(W_R-\mathbb E W_R\), so
\[
 \frac{W_R-\mathbb E W_R}{m_d(B_R)}\longrightarrow 0
 \qquad\text{almost surely}.
\]
The mean bound \(\mathbb E W_R=O_\rho(R^d)\) then gives
\(W_R/m_d(B_R)=O(1)\) almost surely.

For a test function \(f\) supported in \(B_\rho\), write \(T_R(f)\)
for the centered sum with kernel \(Z_{p,q}(f)\) and coefficients
\(c_p=\mu(p)\). Choose a countable family \(\mathcal D_\rho\) that
is dense in the uniform norm among test functions supported in
\(B_\rho\). For each nonzero \(g\in\mathcal D_\rho\), apply the
scalar conclusion to \(Z_{p,q}(g)/\|g\|_\infty\) and multiply
back by \(\|g\|_\infty\). Since the family is countable, there is
a single event of probability one on which
\(T_R(g)/m_d(B_R)\to0\) for every \(g\in\mathcal D_\rho\) and
\(W_R/m_d(B_R)\) is bounded for all sufficiently large \(R\).

Now fix a realization in this event and any test function \(f\)
supported in \(B_\rho\). For \(g\in\mathcal D_\rho\), we have
\[
 |Z_{p,q}(f)-Z_{p,q}(g)|
 \le \|f-g\|_\infty V_{p,q}^{(\rho)}.
\]
Summing this bound with the weights \(|\mu(p)|\,|\mu(q)|\)
controls the difference of the uncentered sums by
\(\|f-g\|_\infty W_R\). Taking expectations gives the corresponding
bound \(\|f-g\|_\infty\mathbb E W_R\) for the difference of their means.
Thus
\[
 \frac{|T_R(f)-T_R(g)|}{m_d(B_R)}
 \le \|f-g\|_\infty\,\frac{W_R+\mathbb E W_R}{m_d(B_R)}.
\]
Since \(T_R(g)/m_d(B_R)\to0\), the triangle inequality gives
\[
 \limsup_{R\to\infty}\frac{|T_R(f)|}{m_d(B_R)}
 \le \|f-g\|_\infty
      \limsup_{R\to\infty}\frac{W_R+\mathbb E W_R}{m_d(B_R)}.
\]
The last limsup is finite and independent of \(g\). By density,
we can make \(\|f-g\|_\infty\) arbitrarily small. This proves
convergence for every test function supported in \(B_\rho\) on
the same event. Finally, intersect these events over positive
integer \(\rho\), since every compact support is contained in
some \(B_\rho\).
\end{proof}

\subsection{Proofs of the auxiliary estimates}

\begin{proof}[Proof of Lemma~\ref{lem:local-pair-bounds}]

Applying Lemma~\ref{lem:local-mark-mass} with
$A_q=\Psi_q$ and $H=B_\rho(x)$ gives
\eqref{eq:mean-local-mark-mass}. Taking $A_q=\Phi_q$
gives the corresponding bound for $\Phi_q$.

Conditioning on \(\mathcal{F}_{(\Phi_p,\Psi_p)}\) fixes \(\Phi_p\).
For \(q\ne p\), the mark \(\Psi_q\) is independent of this
sigma-algebra. Hence
\[
\begin{aligned}
 \mathbb E\!\left[\sum_{q\ne p}|\mu(q)|V_{p,q}^{(\rho)}
                 \,\middle|\,\mathcal{F}_{(\Phi_p,\Psi_p)}\right]
 &=\int\sum_{q\ne p}|\mu(q)|\,\mathbb E \Psi_q(B_\rho(x))\,d\Phi_p(x)\\
 &\le C m_\rho \Phi_p(\R^d)\le C^2m_\rho.
\end{aligned}
\]
The diagonal term is determined by \((\Phi_p,\Psi_p)\) and is at
most \(C^2\|\mu\|_\infty\). Adding it gives the row bound. For the column,
fix \(\Psi_p\) and use the local mass estimate for the \(\Phi_q\).

Under the hypothesis of Theorem~\ref{lem:parent-centering}, the
marks at the other sites are independent of the four measures
at \(p\). Both \(\Phi_p\) and \(\Psi_p\) are fixed under conditioning
on \(\mathcal{F}_{(a_p,b_p,\Phi_p,\Psi_p)}\), so the same calculations
give the row and column bounds under that conditioning.
\end{proof}

\begin{proof}[Proof of Lemma~\ref{lem:window-variance}]
For \(p\in F\), write
\[
 U_p=\sum_{q\in F}|\mu(q)|G_{p,q},
 \qquad
 V_p=\sum_{q\in F}|\mu(q)|G_{q,p}.
\]
Given \(\mathcal{F}_{Y_p}\), the off-diagonal terms in \(U_p\) are
independent, as each depends on a different site \(q\ne p\).
The diagonal term is fixed. Each summand \(W=|\mu(q)|G_{p,q}\)
satisfies \(W^2\le \|\mu\|_\infty L W\). Consequently,
\[
 \operatorname{Var}(U_p\mid\mathcal{F}_{Y_p})
 \le \|\mu\|_\infty L\sum_{q\in F\setminus\{p\}}
          \mathbb E[|\mu(q)|G_{p,q}\mid\mathcal{F}_{Y_p}]
 \le \|\mu\|_\infty L\,\mathbb E[U_p\mid\mathcal{F}_{Y_p}].
\]
The square of the conditional mean is at most \(d_0\) times that
mean. The same argument applies to \(V_p\). Adding the two
second-moment bounds and taking expectations gives
\[
 \mathbb E[U_p^2+V_p^2]\le(\|\mu\|_\infty L+d_0)\mathbb E[U_p+V_p].
\]

In the expansion of \(\operatorname{Var}(S_F)\), the covariance
between \(Z_{p,q}\) and \(Z_{r,s}\) vanishes when the two pairs
have disjoint sites. For the remaining terms we use
\[
 |\operatorname{Cov}(Z_{p,q},Z_{r,s})|
 \le \mathbb E[G_{p,q}G_{r,s}]
       +\mathbb E G_{p,q}\,\mathbb E G_{r,s}.
\]
Bound each
coefficient by its absolute weight and group the resulting
nonnegative terms by a shared index \(i\). The factor contributed
by that index is \(|\mu(i)|^2\).

For example, when \(p=r=i\), the two pairs are \((i,q)\) and
\((i,s)\), so both entries of \(G\) lie in its row \(i\). After
factoring out \(|\mu(i)|^2\), the two parts of the covariance
bound sum over \(q,s\in F\) to
\[
\begin{aligned}
 \sum_{q,s\in F}|\mu(q)|\,|\mu(s)|\,\mathbb E[G_{i,q}G_{i,s}]
 &=\mathbb E[U_i^2],\\
 \sum_{q,s\in F}|\mu(q)|\,|\mu(s)|\,
                \mathbb E[G_{i,q}]\,\mathbb E[G_{i,s}]
 &=(\mathbb E U_i)^2.
\end{aligned}
\]
The first equality comes from expanding \(U_i^2\), and the
second from expanding \((\mathbb E U_i)^2\). By Jensen's
inequality their sum is at most \(2\mathbb E[U_i^2]\).
The table collects this case and the other three, with the
factor \(|\mu(i)|^2\) omitted from each bound.
\[
\begin{array}{c|l}
 \text{Shared index} & \text{Bound}\\ \hline
 p=r=i & \mathbb E U_i^2+(\mathbb E U_i)^2\\
 q=s=i & \mathbb E V_i^2+(\mathbb E V_i)^2\\
 p=s=i & \mathbb E(U_iV_i)+\mathbb E U_i\,\mathbb E V_i\\
 q=r=i & \mathbb E(U_iV_i)+\mathbb E U_i\,\mathbb E V_i.
\end{array}
\]
Every covariance term that can be nonzero is covered by at least
one of these cases for some \(i\in F\). Counting a term more
than once only increases the upper bound, since the entries
are nonnegative.

For a fixed \(i\), the four entries sum to
\(\mathbb E[(U_i+V_i)^2]+(\mathbb E[U_i+V_i])^2\).
By Jensen's inequality and \((u+v)^2\le2(u^2+v^2)\),
\[
\begin{aligned}
 \mathbb E[(U_i+V_i)^2]+(\mathbb E[U_i+V_i])^2
 &\le 2\mathbb E[(U_i+V_i)^2]\le 4\mathbb E[U_i^2+V_i^2].
\end{aligned}
\]
Multiplying this bound by \(|\mu(i)|^2\) and summing over \(i\in F\),
we obtain the first inequality in
\[
\begin{aligned}
 \operatorname{Var}(S_F)
 &\le 4\sum_{i\in F}|\mu(i)|^2\,\mathbb E[U_i^2+V_i^2]\le 4\|\mu\|_\infty(\|\mu\|_\infty L+d_0)\sum_{i\in F}|\mu(i)|\,\mathbb E[U_i+V_i]\\
  &=8\|\mu\|_\infty(\|\mu\|_\infty L+d_0)\,\mathbb E S_F^{+}.
\end{aligned}
\]
The second inequality follows from the second-moment bound for
\(U_i,V_i\) and the estimate \(|\mu(i)|^2\le \|\mu\|_\infty|\mu(i)|\).
For the last equality, we use \(\sum_i|\mu(i)|U_i=S_F^{+}\), which is obtained by summing
over the rows of \(G\),  and \(\sum_i|\mu(i)|V_i=S_F^{+}\), which is obtained by
summing over its columns.

For the full positive sum take \(Z=G\) and \(c_p=|\mu(p)|\).
Since \(S_F^{+}\) and \((S_F^{+})^2\) increase to \(S^{+}\) and \((S^{+})^2\) along a
finite exhaustion, monotone convergence in the finite variance
bound gives the assertion.
\end{proof}

\subsection{Application to random displacements}

Taking Dirac measures as marks, we obtain the two
consequences of Theorems~\ref{soportes}
and~\ref{lem:parent-centering} that will be used
in the study of the random autocorrelation.

\begin{cor}\label{primersoportes}
For any displacement law, the normalized autocorrelations formed
by selecting original sites or final atoms in \(B_R\) differ
vaguely by \(o(1)\) almost surely.
\end{cor}

\begin{proof}
In Theorem~\ref{soportes} take
\(a_p=b_p=\Phi_p=\Psi_p=\delta_{p+\xi_p}\). The observed measures
are both \((\mu_\xi)_R\), and those in the original window are both
\(\sum_{p\in X_R}\mu(p)\delta_{p+\xi_p}\).
\end{proof}

\begin{cor}\label{leyadoc0}
For any displacement law, on one event of probability one, for every
\(f\in C_c(\R^d)\),
    \[
        \lim_{R\to\infty} \frac{1}{m_d(B_R)} \sum_{p,q\in X_R} \mu(p)\overline{\mu(q)} \left[f(p+\xi_p-q-\xi_q) - \e(f(p+\xi_p-q-\xi_q))\right] = 0 \,.
    \]
\end{cor}

\begin{proof}
Apply Theorem~\ref{lem:parent-centering} with
\(a_p=b_p=\Phi_p=\Psi_p=\delta_{p+\xi_p}\). The two original-window
measures coincide with \(\sum_{p\in X_R}\mu(p)\delta_{p+\xi_p}\).
Evaluating their convolution on \(f\) gives the pair sum in
the statement, and centering it gives the limit.
\end{proof}

\section{The random autocorrelation}\label{randomautocorrelation}

The aim of this section is to complete the proof of
Theorem~\ref{autocoraleatoria} and establish
Theorem~\ref{autocoraleatoriaexists}. By Corollary~\ref{primersoportes}, almost surely,
for every $f\in C_c(\R^d)$,
\begin{equation}\label{secanterior}
 \lim_{R\to\infty}\frac{1}{m_d(B_R)}
 \left(
 \bigl((\mu_\xi)_R*\widetilde{(\mu_\xi)_R}\bigr)(f)
 -
 \sum_{p,q\in X_R}\mu(p)\overline{\mu(q)}
 f(p+\xi_p-q-\xi_q)
 \right)=0.
\end{equation}

Together with the centering result of
Corollary~\ref{leyadoc0}, this reduces the study of the
random autocorrelation to the analysis of the expected
pair sums. It is at this stage that the existence of
$\gamma_\mu$ enters.

\begin{proof}[Proof of Theorem~\ref{autocoraleatoria}]
Write $\kappa=\pe_\xi*\widetilde{\pe_\xi}$ and $|B_R|=m_d(B_R)$. A self-pair always contributes $f(0)$. For distinct sites independence gives $\mathbb E f(p-q+\xi_p-\xi_q)=\int f(p-q+u)\,d\kappa(u)$. Thus
\begin{align*}
 \frac1{|B_R|}\sum_{p,q\in X_R}\mu(p)\overline{\mu(q)}\,\mathbb E f(p-q+\xi_p-\xi_q)
 &=\frac1{|B_R|}\sum_{\substack{p,q\in X_R\\p\ne q}}\mu(p)\overline{\mu(q)}
       \int f(p-q+u)\,d\kappa(u)+\frac{\boldsymbol E_R}{|B_R|}f(0)\\
 &=(\gamma_\mu^{(R)}*\kappa)(f)+\frac{\boldsymbol E_R}{|B_R|}\bigl(f(0)-\kappa(f)\bigr).
\end{align*}
The second equality inserts the missing diagonal into the convolution and removes exactly what was inserted. The window comparison and the centering theorem show, on their common event, that
\[
 \gamma_{\mu_\xi}^{(R)}-
 \gamma_\mu^{(R)}*\kappa-
 (\boldsymbol E_R/|B_R|)(\delta_0-\kappa)\longrightarrow0
 \quad\text{vaguely}.
\]
The two energy sums differ only at boundary crossings. Since the atomic weights are bounded, Proposition~\ref{entranysalen} gives
\begin{equation}\label{energylimit}
 |\boldsymbol E_{\xi,R}-\boldsymbol E_R|/|B_R|
 \le \|\mu\|_\infty\,D_R^\xi/|B_R|\longrightarrow0 .
\end{equation}
Finally $\gamma_\mu^{(R)}(f(\cdot+y))\to\gamma_\mu(f(\cdot+y))$ for each $y$. Uniform translation boundedness gives a bound independent of $y$ and $R\ge1$ for fixed $f\in C_c$. Integration against $\kappa$ and dominated convergence finish the proof.
\end{proof}

\begin{proof}[Proof of Theorem~\ref{autocoraleatoriaexists}]

By \eqref{energylimit} 
\begin{equation}\label{existele}
\boldsymbol{E}=\lim_{R\to\infty} \frac{\boldsymbol{E}_{R}}{m_d(B_R)}  \text{ exists } \iff \lim_{R\to\infty} \frac{\boldsymbol{E}_{\xi,R}}{m_d(B_R)}  \text{ exists } \,.
\end{equation}
So, if the limit $\boldsymbol{E}$ exists, by Theorem~\ref{autocoraleatoria}, almost surely $\gamma_{\mu_{\xi}}$ exists, and we have the formula \eqref{eq:random-auto-intro} for the random autocorrelation.

For the converse, assume that $\boldsymbol{E}$ does not exist. Firstly, note that if the common law is not a point mass, there are disjoint open balls with positive probability. Two independent copies fall in the respective balls with positive probability, and on that event their difference is nonzero. Thus for a non-deterministic displacement law we have $\delta_0\ne\pe_\xi*\widetilde{\pe_\xi}$. Choose $f\in C_c(\R^d)$ such that
\[
(\delta_0-\pe_\xi*\widetilde{\pe_\xi})(f)\neq0.
\]
If $\lim_{R\to\infty} \gamma_{\mu_\xi}^{(R)}(f)$ 
existed with positive probability, then Theorem~\ref{autocoraleatoria} would imply the convergence of $\boldsymbol{E}_{\xi,R}/m_d(B_R)$ on an event of positive probability. By \eqref{existele}, this would imply the existence of $\boldsymbol{E}$, a contradiction. Hence, almost surely, $\gamma_{\mu_\xi}$ does not exist.
\end{proof}

\begin{obs}\label{rem:accumulation-segment}
When the diagonal energy density does not exist and the common
displacement law is not a point mass, Theorem~\ref{autocoraleatoriaexists} says that, almost surely, the finite-volume autocorrelations of $\mu_\xi$ have no vague limit along the full family of radii. However, limits along suitable sequences of radii always exist, and all of them can be described explicitly. Indeed, write $\kappa=\pe_\xi*\widetilde{\pe_\xi}$, as in the proof of Theorem~\ref{autocoraleatoria}. By~\eqref{energylimit}, the random energy $\boldsymbol{E}_{\xi,R}$ can be replaced by the deterministic energy $\boldsymbol{E}_R$ in~\eqref{autocorlim}. Hence, almost surely,
\[
    \gamma_{\mu_\xi}^{(R)}-\gamma_\mu*\kappa
    -\frac{\boldsymbol{E}_R}{m_d(B_R)}\,(\delta_0-\kappa)
    \xrightarrow[R\to\infty]{\mathrm{vague}}0.
\]
On the other hand, since $|\mu(p)|^2\leq\|\mu\|_\infty\,|\mu(p)|$ for every $p\in X$, we have that $\boldsymbol{E}_R\leq\|\mu\|_\infty\,|\mu|(B_R)$. So, by~\eqref{eq:TB-poly-growth}, the normalized diagonal energy $\boldsymbol{E}_R/m_d(B_R)$ is bounded for $R\geq1$, and therefore it converges along suitable sequences of radii. If $R_n\to\infty$ is a sequence such that $\boldsymbol{E}_{R_n}/m_d(B_{R_n})$ converges to some $e\geq0$, the above limit shows that, almost surely, the perturbed autocorrelations $\gamma_{\mu_\xi}^{(R_n)}$ converge vaguely to $\gamma_\mu*\kappa+e\,(\delta_0-\kappa)$.

It only remains to determine the possible values of $e$. Let $\boldsymbol{E}_-$ and $\boldsymbol{E}_+$ denote the lower and upper limits of the normalized diagonal energy as $R\to\infty$. The function $R\mapsto\boldsymbol{E}_R$ is a non-decreasing step function, which jumps at those radii $R$ such that $X\cap\partial B_R\neq\varnothing$. Hence, between two consecutive jumps the normalized diagonal energy is continuous and non-increasing, while by Lemma~\ref{lem:TB-shell} the size of its jump at the radius $R\geq1$ is at most
\[
    \|\mu\|_\infty\,\frac{|\mu|(B_R\setminus B_{R-1})}{m_d(B_R)}
    \lesssim\frac1R.
\]
In consequence, given $e\in(\boldsymbol{E}_-,\boldsymbol{E}_+)$, the normalized diagonal energy can cross the level $e$ upwards only by means of a jump. This happens at arbitrarily large radii, and at each of them its distance to $e$ is at most the size of the corresponding jump, which tends to zero. Therefore, every $e\in[\boldsymbol{E}_-,\boldsymbol{E}_+]$ is the limit of the normalized diagonal energy along some sequence of radii, and the vague limits obtained in this way are exactly the measures $\gamma_\mu*\kappa+e\,(\delta_0-\kappa)$ with $e\in[\boldsymbol{E}_-,\boldsymbol{E}_+]$.
\end{obs}
\subsection{Some examples}
Firstly we show that the existence of the autocorrelation does not imply that of the diagonal energy density,
even for a positive translation bounded measure
supported on the union of two uniformly discrete sets.

\begin{exa}\label{ex:diagonal-energy} 
Take $A=\{n\in\mathbb Z:4^j\le |n|<2\cdot4^j\text{ for some integer }j\ge0\}$. Counting complete blocks gives proportions tending to $1/3$ at $R=4^j$ and to $2/3$ at $R=2\cdot4^j$. Define
\[
\mu=\sum_{n\in\mathbb Z\setminus A}\delta_n+
\frac12\sum_{n\in A}
\left(\delta_n+\delta_{n+\varepsilon_n}\right),
\qquad
\varepsilon_n=\frac1{|n|+2} \,.
\]
Each atom of $\delta_{\mathbb Z}$ is either
preserved or replaced by two atoms of total
mass one whose mutual distance tends to zero. 
To compute the autocorrelation, write
$\mu=\sum_{n\in\mathbb Z}\nu_n$, where each
$\nu_n$ is supported in
$[n,n+\varepsilon_n]$ and has total mass one.
For every $f\in C_c(\mathbb R)$ and fixed
$k\in\mathbb Z$, uniform continuity gives
\[
(\nu_n*\widetilde{\nu_{n-k}})(f)
\longrightarrow f(k)
\qquad (|n|\to\infty).
\]
Since $f$ has compact support, only finitely
many values of $k$ contribute to the
convolution sums, uniformly in $R$.
For each fixed $k$,
the number of pairs with $n-m=k$, divided by $2R$,
tends to one. Thus, by Ces\`aro convergence,
the normalized sum of their contributions tends
to $f(k)$, and we obtain
\[
\lim_{R\to\infty}\frac{1}{2R}
\left(
\sum_{|n|\leq R}\nu_n
*
\sum_{|m|\leq R}\widetilde{\nu_m}
\right)(f)
=
\sum_{k\in\mathbb Z}f(k).
\]
Finally, replacing the sums of complete
clusters by $\mu|_{[-R,R]}$ changes the
value of the unnormalized convolution at $f$
by a uniformly bounded amount,
since at most two boundary clusters are
affected and the measures have uniformly
bounded local mass. Hence
$\gamma_\mu=\delta_{\mathbb Z}$.
On the other hand, the diagonal energy
$\boldsymbol E_R$ differs by a uniformly bounded
amount from
\[
\#(\mathbb Z\cap[-R,R])
-\frac12\#(A\cap[-R,R]).
\]
Consequently, the diagonal energy density does not
exist, although the autocorrelation does.
\end{exa}

Theorem~\ref{autocoraleatoriaexists} also covers uniformly discrete sets that are not FLC and whose autocorrelation is not a pure point measure.

\begin{exa}\label{ex:sine-modulation}
 Fix $\alpha\notin\Q$ and $0<c<1/4$, and consider the sine-modulated lattice
\[
    \Lambda=\{x_n:\ n\in\Z\},\qquad x_n=n+c\sin(2\pi n\alpha);
\]
see~\cite[Section~2]{LeeLenzRichardSingStrungaru2020} for this class of modulated structures. For every $n,k\in\Z$ we have that
\[
    x_{n+k}-x_n=k+2c\sin(\pi k\alpha)\cos(2\pi n\alpha+\pi k\alpha).
\]
Taking $k=1$, we see that the gaps between consecutive points are at least $1-2c>1/2$. Hence, $\Lambda$ is uniformly discrete and, as $|x_n-n|\leq c$, it has density one. Moreover, since $\alpha$ is irrational, these gaps take infinitely many different values in the interval $[1/2,3/2]$. Thus, $\Lambda-\Lambda$ is not locally finite, that is, $\Lambda$ is not FLC. To compute the autocorrelation, fix $f\in C_c(\R)$. Since $|x_{n+k}-x_n-k|<1/2$, only finitely many values of $k$ contribute to the autocorrelation sums, and for each of them, the restriction of $x_n$ and $x_{n+k}$ to $[-R,R]$ can be replaced by the condition $|n|\leq R$, changing only a bounded number of terms. Therefore, by Weyl's equidistribution theorem,
\[
    \lim_{R\to\infty}\frac{1}{2R}\sum_{x,y\in\Lambda\cap[-R,R]}f(x-y)
    =\sum_{k\in\Z}\int_0^1f\bigl(k+2c\sin(\pi k\alpha)\cos(2\pi t+\pi k\alpha)\bigr)\,dt.
\]
In other words, $\gamma_{\delta_\Lambda}=\sum_{k\in\Z}\delta_k*\sigma_k$, where $\sigma_k$ denotes the law of $2c\sin(\pi k\alpha)\cos(2\pi t+\pi k\alpha)$, with $t$ uniformly distributed on $[0,1)$. Clearly $\sigma_0=\delta_0$, while for $k\neq0$ the measure $\sigma_k$ is the arcsine law on the interval $(-a_k,a_k)$, where $a_k=2c\,|\sin(\pi k\alpha)|>0$, which is absolutely continuous. So, the autocorrelation exists, and it is not a pure point measure. Since the diagonal energy density of $\delta_\Lambda$ is the density of $\Lambda$, Theorem~\ref{autocoraleatoriaexists} applies for an arbitrary displacement law.
\end{exa}

The window results include displacement laws with infinite first moment.
\begin{exa}
 In dimension one, take the Cauchy density $x\mapsto s/(\pi(s^2+x^2))$, $s>0$. Its characteristic function in our convention is $\varphi_\xi(\lambda)=e^{-2\pi s|\lambda|}$, as follows by integrating this exponential against the inverse Fourier kernel. Thus the law has no atoms and its characteristic function is nowhere zero. For $\mu=\delta_{\Z}$ the observed diffraction converges almost surely to
\[
 \delta_0+\sum_{n\in\Z\setminus\{0\}}e^{-4\pi s|n|}\delta_n
           +(1-e^{-4\pi s|\lambda|})\,m_1.
\]
For the split configuration of Example~\ref{ex:diagonal-energy}, the same law instead destroys the ordinary autocorrelation. The corrected intensities of Corollary~\ref{cor:observable-recovery} recover the common original diffraction $\delta_{\Z}$ in both cases.
\end{exa}

\section{The random diffraction}\label{randomdiffraction}

In this section, we prove Theorems~\ref{thm:uniform-TB} and~\ref{formulasautoanddifrac}. Theorem~\ref{thm:uniform-TB}  establishes uniform translation boundedness of the normalized random autocorrelations, which allows us to strengthen the vague convergence obtained in the previous section to convergence in the space of tempered distributions. This provides the remaining step needed to prove Theorem~\ref{formulasautoanddifrac}: by applying the Fourier transform, we identify the random diffraction and show that it is obtained as a limit from almost every realization. 
For a discussion of the relation between
complex Radon measures and tempered
distributions, see \cite{BaakeStrungaru2023}.

\subsection{Uniform translation boundedness of normalized random
autocorrelations}\label{sec:uniform-TB}

We first establish a uniform bound for the normalized
autocorrelations of the general positive marked measure
introduced in Lemma~\ref{lem:local-mark-mass}.
This will also be used in Section~\ref{generalcase}.

\begin{lem}\label{lem:uniform-mark-TB}
Let $(\alpha_p)_{p\in X}$ be i.i.d.\ random finite positive
measures of mass at most $C$ almost surely, and set
$A_p=\tau_p\alpha_p$. Consider the almost surely Radon
measure
\[
    \rho_A=\sum_{p\in X}|\mu(p)|A_p.
\]
Set, for \(R\geq1\), 
\[ \gamma_{\rho_A}^{(R)}:=\frac{(\rho_A)_R*\widetilde{(\rho_A)_R}}{m_d(B_R)}.
\]
Then, almost surely, the family
$\{\gamma_{\rho_A}^{(R)}:R\ge1\}$ is uniformly translation
bounded.
\end{lem}

The proof of Lemma~\ref{lem:uniform-mark-TB} uses Theorems~\ref{soportes} and~\ref{lem:parent-centering}, already established for every TB purely atomic measure, with locally finite support. Importantly, neither of them requires the deterministic measure to possess autocorrelation. Both results may therefore be applied to $|\mu|$.

\begin{proof}[Proof of Lemma \ref{lem:uniform-mark-TB}]
Theorems~\ref{soportes} and~\ref{lem:parent-centering},
applied to $|\mu|$ with
$a_p=b_p=\Phi_p=\Psi_p=A_p$, give, almost surely,
for every $\varphi\in C_c(\R^d)$,
\[
 \lim_{R\to\infty}\left(
 \gamma_{\rho_A}^{(R)}(\varphi)
 -
 \frac{1}{m_d(B_R)}
 \sum_{p,q\in X_R}|\mu(p)||\mu(q)|
 \mathbb E Z_{p,q}(\varphi)
 \right)=0,
\]
where
\[
 Z_{p,q}(\varphi)
 =\iint\varphi(x-y)\,dA_p(x)\,dA_q(y).
\]
Let $r\ge1$ be such that
$\operatorname{supp}\varphi\subseteq B_r$ and put
$m_r=\sup_x|\mu|(B_r(x))$.
Taking expectations in \eqref{eq:proximity-rows} gives
\[
 \sum_{q\in X_R}|\mu(q)|\,
 \mathbb E|Z_{p,q}(\varphi)|
 \le \|\varphi\|_\infty C^2(m_r+\|\mu\|_\infty).
\]
Consequently,
\[
 \left|
 \frac{1}{m_d(B_R)}
 \sum_{p,q\in X_R}|\mu(p)||\mu(q)|
 \mathbb E Z_{p,q}(\varphi)
 \right| \le
 \|\varphi\|_\infty C^2(m_r+\|\mu\|_\infty)
 \frac{|\mu|(B_R)}{m_d(B_R)}
 \lesssim 1.
\]
This gives a uniform bound on
$|\gamma_{\rho_A}^{(R)}(\varphi)|$ for sufficiently
large $R$. For $1\le R\le R_0$, local finiteness gives
\[
 |\gamma_{\rho_A}^{(R)}(\varphi)|
 \le
 \frac{\|\varphi\|_\infty\rho_A(B_{R_0})^2}{m_d(B_1)}.
\]
Thus
\[
 \sup_{R\ge1}|\gamma_{\rho_A}^{(R)}(\varphi)|<\infty
 \qquad\text{almost surely}.
\]

Choose a nonnegative $u\in C_c(\R^d)$ equal to one
on $B_2$, and set $\varphi:=u*\widetilde u$. For $R\ge1$,
consider
\[
 F_R(x):=\gamma_{\rho_A}^{(R)}(T_x\varphi),
 \qquad T_x\varphi(t):=\varphi(t-x).
\]
To see that $F_R$ is positive definite, take
$x_1,\ldots,x_n\in\R^d$, $c_1,\ldots,c_n\in\C$,
and put $g=\sum_{j=1}^n c_jT_{x_j}u$. Then
\[
 \sum_{i,j=1}^n c_i\overline{c_j}F_R(x_i-x_j)
 =\gamma_{\rho_A}^{(R)}(g*\widetilde g)\ge0,
\]
because $\gamma_{\rho_A}^{(R)}$ is positive definite. Since $\varphi\ge m_d(B_1)\mathbbm1_{B_1}$, we obtain
\[
 m_d(B_1)\gamma_{\rho_A}^{(R)}(B_1(x))
 \le F_R(x)\le F_R(0)
 =\gamma_{\rho_A}^{(R)}(\varphi).
\]
The preceding bound therefore yields
\[
 \sup_{R\ge1}\sup_{x\in\R^d}
 \gamma_{\rho_A}^{(R)}(B_1(x))<\infty
 \qquad\text{almost surely}.
\]
\end{proof}

\begin{proof}[Proof of Theorem~\ref{thm:uniform-TB}]
Apply Lemma~\ref{lem:uniform-mark-TB} with
$A_p=\delta_{p+\xi_p}$, so that $\rho_A=\rho_\xi$.
The conclusion follows from
$|\gamma_{\mu_\xi}^{(R)}|\le\gamma_{\rho_\xi}^{(R)}$.
\end{proof}

A direct consequence of Theorem~\ref{thm:uniform-TB} is the following corollary.

\begin{cor}\label{cor:growth}
Almost surely, there exists a finite constant \(C'(\omega)>0\) such that
\[
    |\gamma_{\mu_\xi}^{(R)}|(B_k)
    \leq C'(\omega)k^d
\]
for every \(R\geq1\) and every \(k\geq1\).
\end{cor}

\subsection{Computing the random diffraction}

As explained above, we now have all the ingredients needed to prove Theorem~\ref{formulasautoanddifrac}.

\begin{proof}[Proof of Theorem~\ref{formulasautoanddifrac}]
If
\[
    \boldsymbol{E}=\lim_{R\to\infty}\frac{1}{m_d(B_R)}\sum_{p\in X_R}|\mu(p)|^2
\]
exists, the right-hand side of
\eqref{eq:random-auto-intro} is translation bounded,
since $\gamma_\mu$ is translation bounded and
$\pe_\xi*\widetilde{\pe_\xi}$ is a finite measure.
Applying Theorem~\ref{autocoraleatoriaexists}, together with Corollary~\ref{cor:growth} and Lemma~\ref{lem:vague-to-distribution}, we have, almost surely
\[
   \lim_{R\to\infty} \gamma_{\mu_\xi}^{(R)}
    =\gamma_\mu*(\pe_\xi*\widetilde{\pe_\xi})
    +\boldsymbol{E}\bigl(\delta_0-\pe_\xi*\widetilde{\pe_\xi}\bigr) = \gamma_{\mu_\xi}
\]
in the distributional sense. In other words, outside a common exceptional null set, the following convergence holds simultaneously for every $f\in \schw$,
\[
\lim_{R\to\infty} \frac{(\mu_{\xi})_R \ast \widetilde{(\mu_{\xi})_R}}{m_d(B_R)} (\widehat{f}) = \bigl(\gamma_\mu*(\pe_\xi*\widetilde{\pe_\xi})
    +\boldsymbol{E}\bigl(\delta_0-\pe_\xi*\widetilde{\pe_\xi}\bigr)\bigr)^{\wedge}(f) = \widehat{\gamma_{\mu_\xi}}(f)\,,
\]
where $(\,\cdot\,)^{\wedge}$ denotes the Fourier transform. Finally
\[
 \bigl(\gamma_\mu*(\pe_\xi*\widetilde{\pe_\xi})
    +\boldsymbol{E}\bigl(\delta_0-\pe_\xi*\widetilde{\pe_\xi}\bigr)\bigr)^{\wedge}(f) = \bigl(|\varphi_\xi|^2\widehat{\gamma_\mu}\bigr)(f)
    +\boldsymbol{E}\int_{\R^d}f(\lambda)\bigl(1-|\varphi_\xi(\lambda)|^2\bigr)\,d\lambda \,.
\]
\end{proof}

\paragraph{Observable corrections and diffraction recovery.}

The correction in Theorem~\ref{autocoraleatoria} involves the random diagonal
energy $\boldsymbol{E}_{\xi,R}$, which need not be determined by the observed
measure $\mu_\xi$ when distinct atoms are displaced to the same location.
Consider instead the \emph{observable diagonal energy}
\[
    \mathcal{E}_{\xi,R}
    :=\sum_{x\in B_R}|\mu_\xi(\{x\})|^2.
\]
Here, the weights of coincident displaced atoms are added before taking the
squared modulus. If $\pe_\xi$ is non-atomic, almost surely no two distinct
atoms collide and hence
$\mathcal{E}_{\xi,R}=\boldsymbol{E}_{\xi,R}$ for every $R>0$.
This observation gives a recovery result that does not require either the
existence of the diagonal energy density $\boldsymbol{E}$ or the original diffraction to be purely atomic.

\begin{cor}\label{cor:observable-recovery}
Suppose that $\gamma_\mu$ exists and that
$\pe_\xi$ is non-atomic. Then, almost surely,
\[
  \widehat{\gamma_{\mu_\xi}^{(R)}}
  -\frac{\mathcal E_{\xi,R}}{m_d(B_R)}
   (1-|\varphi_\xi|^2)m_d
  \xrightarrow[R\to\infty]{\mathcal S'(\R^d)}
  |\varphi_\xi|^2\widehat{\gamma_\mu}.
\]
Consequently, if $\pe_\xi$ is known and $\varphi_\xi$ is nowhere zero, the original diffraction
measure $\widehat{\gamma_\mu}$ can be recovered from a single realization
of $\mu_\xi$.
\end{cor}

\begin{proof}
We have almost surely
$\mathcal E_{\xi,R}=\boldsymbol{E}_{\xi,R}$
for every $R>0$. Moreover, 
\[
\mathcal E_{\xi,R}
\lesssim m_d(B_R), \qquad R\geq1,
\]
almost surely. As before, the corrected autocorrelations
are uniformly translation bounded by
Theorem~\ref{thm:uniform-TB} and the
previous estimate. Theorem
\ref{autocoraleatoria} and Lemma
\ref{lem:vague-to-distribution} then give
convergence in $\mathcal S'(\R^d)$,
exactly as in the proof of Theorem
\ref{formulasautoanddifrac}. Taking Fourier
transforms yields the asserted limit. Write $\eta=|\varphi_\xi|^2\widehat{\gamma_\mu}$ for the limiting measure. If the law is known and $\varphi_\xi$ is nowhere zero, then for each $g\in C_c(\R^d)$ the function $g/|\varphi_\xi|^2$ is continuous with compact support, and
\[
 \widehat{\gamma_\mu}(g)=\eta\bigl(g/|\varphi_\xi|^2\bigr).
\]
Thus $\eta$ determines $\widehat{\gamma_\mu}$.
\end{proof}

\begin{obs*}
    Under the non-atomicity assumption of the preceding corollary, if the diagonal energy density does not exist, Theorem~\ref{autocoraleatoriaexists} implies that $\gamma_{\mu_\xi}$ almost surely does not exist. Nevertheless, the corrected finite-volume diffraction approximants above still converge, allowing recovery of the original diffraction whenever $\varphi_\xi$ is nowhere zero.
\end{obs*}

\section{Independent atomic markings}\label{generalcase}
We return to the independent marks of Theorem~\ref{teo:cluster-main-intro}, and we keep the notation introduced before and in its statement: the marks $\Theta_p$, the bound $C$ on the total absolute descendant weight, the constant $c_\Theta>0$, the energies $\boldsymbol{E}_{\Theta,R}$, the approximants $\gamma_{\mu_\Theta}^{(R)}$ and the measures $M$ and $\kappa$. To include backgrounds whose diagonal energy oscillates, $\boldsymbol{E}_{\Theta,R}$ keeps track of the squared descendant weights in each observed window.
The factor $c_\Theta^{-1}$ in its definition gives $\mathbb E[c_\Theta^{-1}\sum_k|\alpha_k|^2]=1$. Thus, in a window selected by parents, the normalized descendant energy has expectation $\boldsymbol E_R$. In the absence of collisions, $c_\Theta\boldsymbol{E}_{\Theta,R}$ is the sum of the squared observed amplitudes in $B_R$; collisions can require information about the individual descendants.

The mean mark $\mathbb E\Theta$ has total variation at most $C$.
Consequently, $M$ and $\kappa$ have total variation at most $C^2$,
by the bound on the total absolute descendant weight.
These are the two pair measures that enter the marked formula.

For splitting, even a deterministic separation within a mark contributes a phase difference to $\widehat\Theta$. For example, let $B$ be Bernoulli with mean $\theta$. With $\Theta=\delta_0+B\delta_v$, $v\ne0$, one has
\[
 \mathbb E\widehat\Theta(\lambda)=1+\theta e^{-2\pi i\langle v,\lambda\rangle},
 \qquad \operatorname{Var}(\widehat\Theta(\lambda))=\theta(1-\theta).
\]
The displacement $v$ changes the coherent modulation, while the diffuse term measures the fluctuation in the second descendant's occupation.

The general estimates of Sections~\ref{randomperturbTB}
and~\ref{averaging} apply directly to the marked setting
by counting crossings with their absolute descendant
weights. Two additional points require attention:
the fluctuations of the diagonal energy involve an
entire mark, and descendants of the same parent site
need not be independent. Throughout this section,
the marks at distinct parent sites are independent.

Consider the positive Borel measure
\[
\rho_\Theta
:=
\sum_{p\in X}|\mu(p)|
\sum_{k=1}^{N(p)}|\alpha_{p,k}|\,\delta_{p+\xi_{p,k}}\,.
\]

\begin{prop}\label{radongeneral}
    Almost surely, both $\rho_\Theta$ and $\mu_{\Theta}$ are Radon measures and satisfy $|\mu_{\Theta}|\leq \rho_\Theta$.
\end{prop}

\begin{proof}
Apply Lemma~\ref{lem:local-mark-mass} with
\[
    A_p=\sum_{k=1}^{N(p)}
         |\alpha_{p,k}|\delta_{p+\xi_{p,k}},
    \qquad a_p=\tau_p\Theta_p.
\]
By the bound on the total absolute descendant weight,
$A_p(\R^d)\le C$, and clearly $|a_p|\le A_p$.
The corresponding measures are
$\rho_A=\rho_\Theta$ and $\nu=\mu_\Theta$.
\end{proof}

The marked measures satisfy the same local growth estimate as
in Theorem~\ref{emeR}.

\begin{teo}\label{emeRgeneral}
Let \(D>d\,C\|\mu\|_\infty\). Almost surely,
\(\rho_\Theta(Q+q)\le D\,\ell(\|q\|)\) for all
\(q\in\Z^d\) of sufficiently large norm. In particular, \(\rho_\Theta\)
and \(\mu_\Theta\) are quasi-TB.
\end{teo}

\begin{proof}
For each \(q\), \(\rho_\Theta(Q+q)\) is a sum of independent site
contributions
\[
 |\mu(p)|\sum_k|\alpha_{p,k}|\mathbbm1_{Q+q}(p+\xi_{p,k}),
\]
bounded by \(C\|\mu\|_\infty\). By Lemma~\ref{lem:local-mark-mass},
\[
 \mathbb E\rho_\Theta(Q+q)
 \le C\sup_{x\in\R^d}|\mu|(Q+x)<\infty.
\]
The bound is uniform in \(q\). Lemma~\ref{lem:lattice-growth},
applied with \(b=C\|\mu\|_\infty\) and \(a=D/b>d\),
therefore gives the asserted bound, and
\(|\mu_\Theta|\le\rho_\Theta\) transfers the quasi-TB property.
\end{proof}

The discrepancy between the two windows is measured by the total absolute weight of crossing descendants. Put \[
 D_R^\Theta := \sum_{p\in X} |\mu(p)| \sum_{k=1}^{N(p)} |\a_{p,k}| |\mathbbm{1}_{\{p+\xi_{p,k}\in B_R\}} - \mathbbm{1}_{\{p\in B_R\}}|  \,.
\]
\begin{prop}\label{entranysalengeneral}
Almost surely, \(D_R^\Theta=o(R^d)\) as \(R\to\infty\).
\end{prop}

\begin{proof}
Apply Lemma~\ref{coronas} with
\(A_p=\sum_k|\alpha_{p,k}|\delta_{p+\xi_{p,k}}\).
These are translated positive marks of mass at most \(C\), and
\(D_R^\Theta=D_R^A\).
\end{proof}

\begin{cor}
 Almost surely, \(\rho_\Theta(B_R)\lesssim R^d\) as $R\to\infty$.
\end{cor}
\begin{proof}
The descendants of parents in $B_R$ contribute at most $C|\mu|(B_R)$ to $\rho_\Theta(B_R)$. Every remaining descendant in the ball is counted in $D_R^\Theta$, so Proposition~\ref{entranysalengeneral} and the deterministic volume bound give the result.
\end{proof}

Now
\[
 \bigl((\mu_{\Theta})_R \ast \widetilde{(\mu_{\Theta})_R}\bigr)(f) = \sum_{p,q\in X} \mu(p)\overline{\mu(q)}\sum_{k=1}^{N(p)}\sum_{l=1}^{N(q)} \a_{p,k}\overline{\a_{q,l}} \mathbbm{1}_{\{p+\xi_{p,k}\in B_R\}} \mathbbm{1}_{\{q+\xi_{q,l}\in B_R\}} f(p+\xi_{p,k}-q-\xi_{q,l})  \,.
\]

Let 
\[
 Z_{p,q}(f) := \sum_{k=1}^{N(p)}\sum_{l=1}^{N(q)} \a_{p,k}\overline{\a_{q,l}} f(p+\xi_{p,k}-q-\xi_{q,l}) \,.
\]
We can now compare the observed window with the window indexed by the parent sites.
\begin{cor}\label{soportesgeneral}
Almost surely, for every $f\in C_c(\R^d)$,
\begin{equation}\label{ecsoportesgeneral}
 \frac{1}{m_d(B_R)}\left(\bigl((\mu_{\Theta})_R * \widetilde{(\mu_{\Theta})_R}\bigr)(f)
       - \sum_{p,q\in X_R}\mu(p)\overline{\mu(q)}Z_{p,q}(f)\right)
 \longrightarrow 0.
\end{equation}
\end{cor}
\begin{proof}
In Theorem~\ref{soportes} take $a_p=b_p=\tau_p\Theta_p$ and
$\Phi_p=\Psi_p=\sum_k|\alpha_{p,k}|\,\delta_{p+\xi_{p,k}}$. The domination also holds if descendants coincide, and the bound on the total absolute weight gives \(\Phi_p(\R^d)\le C\). The pairs \((\Phi_p,\Psi_p)\) are independent across sites. For this choice, the parent-window convolution is exactly the second sum in \eqref{ecsoportesgeneral} and the observed-window convolution is the first.
\end{proof}

After this replacement, centering the pair sum gives the following strong law.

\begin{teo}\label{leyadoc0general}
 Under the bound on the total absolute descendant weight alone, almost surely, for each $f\in C_c(\R^d)$
    \[
        \lim_{R\to\infty} \frac{1}{m_d(B_R)}  \sum_{p,q\in X_R}  \mu(p)\overline{\mu(q)}\left(Z_{p,q}(f)-\e(Z_{p,q}(f))\right)= 0 \,.
    \]
\end{teo}

\begin{proof}[Proof of Theorem~\ref{leyadoc0general}]
Take \(a_p=b_p=\tau_p\Theta_p\) and
\(\Phi_p=\Psi_p=\sum_k|\alpha_{p,k}|\delta_{p+\xi_{p,k}}\)
in Theorem~\ref{lem:parent-centering}. Each quadruple
\((a_p,b_p,\Phi_p,\Psi_p)\) is a function of the descendant
data at \(p\), so the quadruples are independent across sites. The total
absolute weight bounds the mass of each positive mark by \(C\). The resulting kernel is precisely \(Z_{p,q}(f)\).
\end{proof}

At a marked site the sum of squared descendant weights fluctuates even when every descendant stays in the window. The following lemma separates these fluctuations from the boundary error.

\begin{lem}\label{energyclustergeneral}
    Almost surely,
\[
 \lim_{R\to\infty}
 \frac{\boldsymbol E_{\Theta,R}-\boldsymbol E_R}
 {m_d(B_R)}
 =0.
\]
\end{lem}

\begin{proof}
Put $\zeta_p=\sum_k|\alpha_{p,k}|^2$. These variables are independent, bounded by $C^2$ and have mean $c_\Theta>0$. Selecting descendants by their parent sites first gives
\[
 c_\Theta(\boldsymbol E_{\Theta,R}-\boldsymbol E_R)
 =\sum_{p\in X_R}|\mu(p)|^2(\zeta_p-c_\Theta)+H_R,
 \qquad |H_R|\le C\|\mu\|_\infty\,D_R^\Theta .
\]
To see this, write $I_R(x)=\mathbbm1_{B_R}(x)$. The error from replacing the final location by the original site is exactly
\[
 H_R=\sum_{p\in X}|\mu(p)|^2\sum_k|\alpha_{p,k}|^2
        \bigl(I_R(p+\xi_{p,k})-I_R(p)\bigr).
\]
Taking absolute values and using $|\alpha_{p,k}|^2\le C|\alpha_{p,k}|$ gives the stated bound. Proposition~\ref{entranysalengeneral} gives $H_R=o(R^d)$ almost surely.

Write \(T_R=\sum_{p\in X_R}|\mu(p)|^2(\zeta_p-c_\Theta)\).
The terms are centered and independent, and
\(0\le \zeta_p,c_\Theta\le C^2\), so
\[
 \operatorname{Var}(T_R)
 \le C^4\sum_{p\in X_R}|\mu(p)|^4
 \le C^4 \|\mu\|_\infty^3\,|\mu|(B_R)=O(R^d).
\]
Lemma~\ref{lem:quadratic-summability} gives \(T_{r_n}=o(r_n^d)\)
almost surely for \(r_n=n^2\). For \(r_n\le R\le s_n=r_{n+1}\),
\[
 |T_R-T_{r_n}|
 \le C^2 \|\mu\|_\infty\,|\mu|(B_{s_n}\setminus B_{r_n})=o(r_n^d).
\]
Together with the estimate on \(H_R\), this proves the lemma.
\end{proof}

\subsection*{Identification of autocorrelation and diffraction}

\begin{proof}[Proof of \eqref{firstpart} and \eqref{secondpart}]
Independence at distinct sites gives $\mathbb E Z_{p,q}(f)=\kappa(f(p-q+\cdot))$ for $p\ne q$. At $p=q$ both descendants belong to the same mark, giving $M(f)$. Thus the expected pair sum divided by $m_d(B_R)$ equals
\[
 (\gamma_\mu^{(R)}*\kappa)(f)+
 \frac{\boldsymbol E_R}{m_d(B_R)}(M-\kappa)(f).
\]
Corollary~\ref{soportesgeneral} and Theorem~\ref{leyadoc0general} show that this differs from the observed autocorrelation at $f$ by $o(1)$, simultaneously for all test functions. Uniform translation boundedness permits dominated convergence against the finite measure $|\kappa|$, exactly as in the proof of Theorem~\ref{autocoraleatoria}. Lemma~\ref{energyclustergeneral} then gives the corrected limit, and adding back the energy term gives \eqref{secondpart} whenever its density exists.
\end{proof}

\begin{proof}[Proof of \eqref{thirdpart} and \eqref{eq:computable-cluster-diffraction}]

Apply Lemma~\ref{lem:uniform-mark-TB} with
\[
 A_p=\sum_{k=1}^{N(p)}
 |\alpha_{p,k}|\delta_{p+\xi_{p,k}},
\]
so that $\rho_A=\rho_\Theta$. Since
$|\gamma_{\mu_\Theta}^{(R)}|\le\gamma_{\rho_\Theta}^{(R)}$,
the family $\{\gamma_{\mu_\Theta}^{(R)}:R\ge1\}$ is
almost surely uniformly translation bounded. In particular,
it satisfies the growth estimate in
Lemma~\ref{lem:vague-to-distribution}. When $\boldsymbol E$ exists, the limit in \eqref{secondpart} is translation bounded because $\kappa$ and $M$ are finite. The convergence therefore holds against Schwartz functions. Applying it to $\widehat f$ gives convergence of the observed scattering at $f$. Finally Fubini, justified by $\|\Theta\|\le C$, gives
\[
 \widehat M(\lambda)=\mathbb E|\widehat\Theta(\lambda)|^2,\qquad
 \widehat \kappa(\lambda)=|\mathbb E\widehat\Theta(\lambda)|^2
\]
as Fourier--Stieltjes transforms. Substitution into \eqref{eq:FT-conv-mult} and \eqref{secondpart} proves both formulas.
\end{proof}

\paragraph{The exceptional mark laws.}
There is a marked analogue of the deterministic exception in Theorem~\ref{autocoraleatoriaexists}. Here it is the mark as a measure that matters. Relabeling the same descendants at random, for instance, need not introduce any randomness in this measure.

\begin{cor}\label{cor:marked-obstruction}
Under the hypotheses of Theorem~\ref{teo:cluster-main-intro}, suppose the diagonal energy density does not exist. If the law of $\Theta$ as a finite complex measure is not concentrated on one measure, then the marked autocorrelation almost surely does not exist. If $\Theta$ is almost surely a fixed measure $\theta$, the marked autocorrelation exists and equals $\gamma_\mu*(\theta*\widetilde\theta)$.
\end{cor}
\begin{proof}
Write $C_\Theta=M-\kappa$. Lemma~\ref{energyclustergeneral} and \eqref{firstpart} give, on one full-measure event and for every $g\in C_c(\R^d)$,
\[
 \gamma_{\mu_\Theta}^{(R)}(g)
 =(\gamma_\mu*\kappa)(g)+\frac{\boldsymbol E_R}{m_d(B_R)}\,C_\Theta(g)+o_g(1).
\]
If $C_\Theta\ne0$, choose $g$ with $C_\Theta(g)\ne0$. The right-hand side has two distinct subsequential limits, because the normalized energy is bounded and fails to converge. If $C_\Theta=0$, the same identity directly gives convergence.

It remains to identify which marks have $C_\Theta=0$. For each $\lambda$ the Fourier--Stieltjes transform satisfies
\[
 \widehat{C_\Theta}(\lambda)
 =\mathbb E\,|\widehat\Theta(\lambda)-\mathbb E\widehat\Theta(\lambda)|^2 .
\]
If $C_\Theta=0$, this variance vanishes at all rational $\lambda\in\mathbb Q^d$. Taking the countable intersection gives one event on which $\widehat\Theta$ agrees with $\widehat{\mathbb E\Theta}$ on $\mathbb Q^d$. Both are continuous transforms of finite measures, so they agree everywhere. Uniqueness of the Fourier transform gives $\Theta=\mathbb E\Theta$ on that event. Conversely, if $\Theta=\theta$ almost surely, then $M=\kappa=\theta*\widetilde\theta$ by definition. This proves both assertions.
\end{proof}

\paragraph{What can be read from a marked realization.}
For a known mark law, the quantities $c_\Theta$, $\kappa$ and $M$ are deterministic. Assume that almost surely no two nonzero descendants occupy the same position, whether from the same parent or from different parents. Absolute local summability then gives, simultaneously for all $R$,
\[
 \boldsymbol E_{\Theta,R}
   =c_\Theta^{-1}\sum_{x\in B_R}|\mu_\Theta(\{x\})|^2 .
\]
This quantity is $O(m_d(B_R))$ almost surely, by Lemma~\ref{energyclustergeneral} and boundedness of $\boldsymbol E_R/m_d(B_R)$ at large radii, and by local finiteness on bounded ranges of radii. The preceding uniform translation bound for $\gamma_{\mu_\Theta}^{(R)}$ used only the estimates for the positive majorant, so it also holds when $\boldsymbol E$ fails to exist. Subtracting the bounded multiples of the finite measure $M-\kappa$ preserves this uniform bound. Lemma~\ref{lem:vague-to-distribution} applied to \eqref{firstpart} therefore gives, after Fourier transformation,
\[
 \widehat{\gamma_{\mu_\Theta}^{(R)}}-
 \frac{\boldsymbol E_{\Theta,R}}{m_d(B_R)}
 \operatorname{Var}(\widehat\Theta)\,m_d
 \ \longrightarrow\ |\fii_\Theta|^2\,\widehat{\gamma_\mu}
 \quad\text{in }\mathcal S' .
\]
If $\fii_\Theta$ is nowhere zero, its reciprocal is bounded on each compact set. The limiting measure thus determines the original diffraction by local division, just as for displacements. Atomic marks need not satisfy the collision assumption. For instance, two coincident descendants with weights $1$ and $-1$ contribute two to the sum of squared descendant weights per unit parent amplitude but disappear from the observed measure. This explains the role of the assumption in this observable formula.

\subsection*{Spectral hyperuniformity.}
Recall that a measure $\nu$ whose autocorrelation exists is called \emph{spectrally hyperuniform} if
\[
    \widehat{\gamma_\nu}(B_r\setminus\{0\})=o(m_d(B_r))
    \qquad(r\to0^+);
\]
see~\cite{BaakeGrimmHU,BjorklundHartnickHU,BjorklundBylehn2026,LachiezeRey2025} for more information about this notion, and~\cite{DereudreFlimmelHuesmannLeble} for dependent perturbations of lattices, where hyperuniformity may fail to be preserved. The diffraction formula~\eqref{thirdpart} has the following consequence.

\begin{cor}\label{cor:spectral-hu}
Under the hypotheses of Theorem~\ref{teo:cluster-main-intro},
assume that the diagonal energy density $\boldsymbol{E}$ exists.
If $\mu$ is spectrally hyperuniform, then, almost surely,
\[
    \lim_{r\to0^+}
    \frac{\widehat{\gamma_{\mu_\Theta}}(B_r\setminus\{0\})}{m_d(B_r)}
    =\boldsymbol{E}\,\operatorname{Var}\Bigl(\,\sum_{k=1}^{N}\alpha_k\Bigr).
\]
For pure displacements:
$\mu_\xi$ is almost surely spectrally hyperuniform
if and only if $\mu$ is spectrally hyperuniform.
More generally, if $\boldsymbol{E}>0$ and
$\mu$ is spectrally hyperuniform, the measure
$\mu_\Theta$ is almost surely spectrally hyperuniform
if and only if the total weight
$\sum_{k=1}^{N}\alpha_k$ of the mark
is almost surely constant.
\end{cor}

\begin{proof}
By~\eqref{thirdpart}, almost surely, for every $r>0$ we have that
\[
    \widehat{\gamma_{\mu_\Theta}}(B_r\setminus\{0\})
    =\int_{B_r\setminus\{0\}}|\fii_\Theta|^2\,d\widehat{\gamma_\mu}
    +\boldsymbol{E}\int_{B_r}\operatorname{Var}(\widehat\Theta)(\lambda)\,d\lambda.
\]
By~\eqref{eq:mark-bound}, we have that $|\widehat\Theta(\lambda)|\leq C$ almost surely, for every $\lambda\in\R^d$. On the one hand, this implies that $|\fii_\Theta|\leq C$. So, the first term of the right-hand side is at most $C^2\,\widehat{\gamma_\mu}(B_r\setminus\{0\})$, which is $o(m_d(B_r))$ because $\mu$ is spectrally hyperuniform. On the other hand, by the dominated convergence theorem, the function $\operatorname{Var}(\widehat\Theta)$ is continuous, and its value at the origin is the variance of $\widehat\Theta(0)=\sum_{k=1}^{N}\alpha_k$. This proves the formula of the statement. The forward implication for $\mu_\xi$ follows by taking
$\Theta=\delta_\xi$, whose total weight is equal to one.
Conversely, by~\eqref{eq:random-diff-intro} and positivity,
\[
\widehat{\gamma_{\mu_\xi}}(B_r\setminus\{0\})
\geq
\left(\inf_{\lambda\in B_r}|\varphi_\xi(\lambda)|^2\right)
\widehat{\gamma_\mu}(B_r\setminus\{0\}).
\]
Since $\varphi_\xi$ is continuous and
$\varphi_\xi(0)=1$, the infimum tends to one
as $r\to0^+$. Hence spectral hyperuniformity
of $\mu_\xi$ implies that of $\mu$.
Finally, when $\boldsymbol{E}>0$ and $\mu$ is
spectrally hyperuniform, the last assertion
follows because the variance of a random variable
is zero if and only if it is almost surely constant.
\end{proof}
For instance, if $\boldsymbol{E}>0$ and $\mu$ is
spectrally hyperuniform, spectral hyperuniformity is destroyed by the random occupation $\Theta=\alpha\,\delta_0$, where $\alpha$ is a Bernoulli random variable with mean $\theta\in(0,1)$, whereas it is preserved when each atom is split into several descendants whose weights add up to a fixed constant. This is the case, for instance, when the weights are random convex combinations.

\subsection*{General finite marks}

The atomic structure of the marks is not essential for the
identification of autocorrelation and diffraction. Assume that
$\gamma_\mu$ exists and let $(\Theta_p)_{p\in X}$ be i.i.d.\ random
finite complex measures, not necessarily atomic, with the same law as $\Theta$, satisfying
$|\Theta|(\R^d)\leq C$ almost surely. Note that \eqref{eq:mark-bound} implies this bound,
but the converse need not hold when descendants
of the same parent occupy the same position. Define
$\mu_\Theta=\sum_{p\in X}\mu(p)\tau_p\Theta_p$, and retain the
notation $M$, $\kappa$ and $\fii_\Theta=\mathbb E\widehat\Theta$.

Set $A_p=\tau_p|\Theta_p|$. The total variation is measurable
as a function of the mark, since its value on any Borel set can
be expressed as a supremum over finite partitions induced
by a countable generating algebra. Lemma~\ref{lem:local-mark-mass}
then gives local finiteness and domination.
Theorems~\ref{soportes} and~\ref{lem:parent-centering} apply with
$a_p=b_p=\tau_p\Theta_p$ and $\Phi_p=\Psi_p=A_p$.
The computation of the expected pair sum in the proof of
\eqref{firstpart} remains unchanged, and the same dominated
convergence argument gives, almost surely,
\[
 \gamma_{\mu_\Theta}^{(R)}-
 \frac{\boldsymbol E_R}{m_d(B_R)}(M-\kappa)
 \xrightarrow[R\to\infty]{\mathrm{vague}}
 \gamma_\mu*\kappa.
\]
By Lemma~\ref{lem:uniform-mark-TB} and the domination
$|\gamma_{\mu_\Theta}^{(R)}|\leq\gamma_{\rho_A}^{(R)}$,
the corrected approximants are uniformly translation bounded,
since $\boldsymbol E_R/m_d(B_R)$ is bounded. Thus the convergence
also holds in $\mathcal S'(\R^d)$ by
Lemma~\ref{lem:vague-to-distribution}.

In particular, \eqref{secondpart}, \eqref{thirdpart} and
\eqref{eq:computable-cluster-diffraction} hold whenever
$\boldsymbol E$ exists, and the conclusions of
Corollary~\ref{cor:marked-obstruction} extend to this setting.
The marked assertions of Corollary~\ref{cor:spectral-hu}
also remain valid, with $\Theta(\R^d)$ in place of
$\sum_k\alpha_k$. The atomic structure of the marks is,
however, relevant to the observable correction based on
squared descendant weights, which does not automatically
extend to this setting.

\subsection*{Use of artificial intelligence}

During the preparation of this work the authors used
AI-tools as follows.

\emph{Mathematical content.} The present version of the proof of
Theorem~\ref{thm:uniform-TB}, which upgrades vague convergence
to convergence in $\mathcal{S}'(\R^d)$, was produced with the help of the AI-tools; 
earlier drafts by the authors reached the same
conclusions by a longer argument. The tool also contributed to
the moment-free centering arguments.

\emph{Writing.} The tool was used to improve the exposition
and the English.

All AI-generated material was checked in detail, and corrected
or rewritten where necessary, by the authors, who also
verified the cited literature independently. The authors take
full responsibility for the content of this article. The AI-tools
are not authors.

\end{document}